\documentclass[a4paper,leqno,10pt]{amsart}

\usepackage{epsf,graphicx,epsfig}
\usepackage{amscd}
\usepackage{amsmath,latexsym,amssymb,amsthm}
\usepackage[nospace,noadjust]{cite}
\usepackage{textcomp}
\usepackage{setspace,cite}
\usepackage{lscape,fancyhdr,fancybox}
\usepackage{stmaryrd}
\usepackage[all,cmtip]{xy}
\RequirePackage{tikz-cd}
\RequirePackage{amssymb}
\usetikzlibrary{calc}
\usetikzlibrary{decorations.pathmorphing}

\tikzset{curve/.style={settings={#1},to path={(\tikztostart)
			.. controls ($(\tikztostart)!\pv{pos}!(\tikztotarget)!\pv{height}!270:(\tikztotarget)$)
			and ($(\tikztostart)!1-\pv{pos}!(\tikztotarget)!\pv{height}!270:(\tikztotarget)$)
			.. (\tikztotarget)\tikztonodes}},
	settings/.code={\tikzset{quiver/.cd,#1}
		\def\pv##1{\pgfkeysvalueof{/tikz/quiver/##1}}},
	quiver/.cd,pos/.initial=0.35,height/.initial=0}

\tikzset{tail reversed/.code={\pgfsetarrowsstart{tikzcd to}}}
\tikzset{2tail/.code={\pgfsetarrowsstart{Implies[reversed]}}}
\tikzset{2tail reversed/.code={\pgfsetarrowsstart{Implies}}}
\tikzset{no body/.style={/tikz/dash pattern=on 0 off 1mm}}

\usepackage{graphicx}

\usepackage{color}
\usepackage{url}
\usepackage{enumerate}
\usepackage[mathscr]{euscript}

\newcounter{cnt1}
\newcounter{cnt2}
\newcounter{cnt3}
\newcommand{\blr}{\begin{list}{$($\roman{cnt1}$)$} {\usecounter{cnt1}
			\setlength{\topsep}{0pt} \setlength{\itemsep}{0pt}}}
	\newcommand{\bla}{\begin{list}{$($\alph{cnt2}$)$} {\usecounter{cnt2}
				\setlength{\topsep}{0pt} \setlength{\itemsep}{0pt}}}
		\newcommand{\bln}{\begin{list}{$($\arabic{cnt3}$)$} {\usecounter{cnt3}
					\setlength{\topsep}{0pt} \setlength{\itemsep}{0pt}}}
			\newcommand{\el}{\end{list}}
		
		\newtheorem{Thm}{Theorem}[section]
		
		\newtheorem{Prop}[Thm]{Proposition}
		\newtheorem{Def}[Thm]{Definition}
		\newtheorem{Exm}[Thm]{Example}
		\newtheorem{Rem}[Thm]{Remark}
		\newtheorem{Cor}[Thm]{Corollary}
		
        \newtheorem{Nota}[Thm]{Notation}

		\title{}
		\author{}
		\date{}
		\usepackage{amssymb}

		\usepackage{hyperref}
		\hypersetup{
			colorlinks,
			citecolor=blue,
			filecolor=black,
			linkcolor=blue,
			urlcolor=black
		}
\begin{document}
\title{Hermitian Lie algebroids over analytic spaces}
\author{Abhishek Sarkar}
\begin{abstract}
We explore complex Riemannian geometry and Hermitian metrics on complex algebraic varieties and analytic spaces, respectively. In particular, we introduce Hermitian metrics on holomorphic Lie algebroids and examine the associated characteristic foliation with its canonically induced inner product. Furthermore, we study hypercohomologies related to the leaf space, leaves, and certain invariant subspaces arising from the characteristic foliation of a holomorphic Lie algebroid over a Hermitian manifold. Finally, we extends the concept of equivariant de Rham cohomology to the setting of Hermitian Lie algebroids.
\end{abstract}
\maketitle
\section{Introduction}
Riemannian manifolds are a classical object in differential geometry, defined as smooth manifolds equipped with a Riemannian metric on their tangent bundles, assigning a smoothly varying positive-definite inner product on each tangent space \cite{SM}. In complex geometry, analogous notions give rise to Hermitian manifolds \cite{GH} and holomorphic Riemannian manifolds \cite{Le,DZ}. More recently, M.~Boucetta introduced the notion of Riemannian metrics on smooth Lie algebroids \cite{MB}. These are Lie algebroids over a smooth manifold whose underlying smooth vector bundles are equipped with a Riemannian metric, and which satisfy a canonical compatibility condition with the induced Levi-Civita connection. The study of generalized infinitesimal symmetries led to the notion of (smooth) Lie algebroids in the context of Lie groupoids, which generalize the tangent bundle and Lie algebras (see \cite{KM,RF}). More generally, their algebraic counterparts are known as Lie-Rinehart algebras \cite{GR,JH}. Subsequently, Lie algebroids were developed in holomorphic settings over complex manifolds \cite{SC,LG-PX,PT} and, more generally, over analytic spaces by B.~Pym \cite{BP}. In algebro-geometric contexts, Lie algebroids provide a unified framework for the study of geometric structures such as Poisson analytic spaces, meromorphic connections, singular foliations, and free Lie algebroids \cite{MK, BP, JV, AA, TS}.

In this work, we extend the metric theories from complex manifolds to analytic spaces and complex algebraic varieties. More precisely, we develop foundational aspects of complex Riemannian geometry for complex algebraic varieties and of Hermitian geometry on analytic spaces, building on ideas from complex Riemannian geometry \cite{Le,DZ}, Riemannian Rinehart spaces \cite{VP}, and classical Hermitian geometry \cite{GH}. We study an analytic subspace $Y$ of a Hermitian manifold $X$ and establish conditions under which $Y$ inherits a natural metric from $X$ (see Theorem~\ref{induced inner product}). Our approach relies on \emph{coherent sheaves of Lie-Rinehart algebras} (i.e., Lie algebroids), in particular the tangent sheaves $\mathcal{T}_X$ and $\mathcal{T}_Y$ of $X$ and $Y$, respectively, together with the sheaf of logarithmic derivations (or log tangent sheaf) $\mathcal{T}_X(-\log Y)$ for $Y \subset X$ \cite{BP,AA}. Motivated by the notion of the Levi-Civita connection on a Riemannian manifold, which is the unique torsion-free affine connection preserving the metric \cite{SM}, we construct an analogue of the Levi-Civita connection in the Hermitian setting, yielding a covariant connection on the tangent sheaf $\mathcal{T}_Y$ of $Y$ that is compatible with the induced metric (see Theorem~\ref{connection on $Y$}). Locally, this construction involves the standard Hermitian metric and the corresponding Levi-Civita connection on complex manifolds via the usual covariant derivative. We globalize these locally defined inner products and connections using a smooth partition of unity subordinate to an atlas of $X$. Altogether, this framework provides an analogue of Riemannian geometry on smooth manifolds in the broader context of singular spaces, formulated within an algebro-geometric setting.

The notion of Atiyah algebroid of a vector bundle, with the canonical Lie algebroid structure, has many applications in literature \cite{PT2, PT}. More generally, for a (quasi-)coherent $\mathcal{O}_X$-module $\mathcal{E}$, we have the notion of \emph{Atiyah algebroid} of $\mathcal{E}$, denoted by $\mathcal{A}t(\mathcal{E})$, consists of sheaf of first order differential operators on $\mathcal{E}$~ \cite{MK,BRT,Ab}. In order to study \emph{logarithmic connections} for a smooth complex divisor $Y \subset X$, P. Tortella considered modules over the Atiyah algebroid on the normal bundle of $Y$ in $X$ \cite{PT}. Here, we derive a relationship between the Atiyah algebroid $\mathcal{A}t(\mathcal{I})$ over a principal ideal sheaf $\mathcal{I} \subset \mathcal{O}_X$ \cite{Ab} and sheaf of logarithmic derivations $\mathcal{T}_X(-log ~Y)$ along the associated principal divisor $(Y:=V(\mathcal{I}),~\mathcal{O}_Y:= \mathcal{O}_X/{\mathcal{I}})$ of a Hermitian manifold $X$ (see Theorem \ref{log tangent-Atiyah sheaf}). Here, $Y \subset X$ is obtained by the vanishing locus for the ideal sheaf $\mathcal{I}$.

The notion of Riemannian and Hermitian metrics on real smooth and holomorphic vector bundles, respectively, is well established in the literature; see, for instance, \cite{SM, GH}. Moreover, Hermitian metrics may also be defined on (almost) complex vector bundles. In a related direction, the notion of Hermitian metrics on almost complex smooth Lie algebroids \cite{IP} arises as an analogue of the concept of Riemannian metrics on smooth Lie algebroids \cite{MB}. Although the focus of the present article is different, since we work entirely in the holomorphic setting, dealing either with holomorphic vector bundles or with coherent $\mathcal{O}_X$-modules over a complex manifold $X$, these ideas naturally motivate the definition of Hermitian metrics on holomorphic Lie algebroids. More generally, this framework extends to Lie algebroids over analytic spaces in the algebro-geometric language, which we refer to as \emph{Hermitian Lie algebroids} (see Definition~\ref{definition of hermitian Lie algebroids}).

The notion of characteristic (or orbit) foliation of a smooth Lie algebroid $\mathfrak{a}\colon L \rightarrow TX$ is defined by the image of the anchor map $\mathfrak{a}(L)$, or $\operatorname{Im}(\mathfrak{a})$, which describes a Stefan-Sussmann (or generalized involutive) distribution \cite{RF,MB}. M.~Boucetta studied the induced Riemannian manifold structure on a leaf associated with the characteristic foliation of a Lie algebroid equipped with a Riemannian metric \cite{MB}. We start with a holomorphic Lie algebroid and study the associated \emph{characteristic (analytic) foliation}. In particular, we consider leaves (or orbits) and $\mathcal{L}$-invariant subspaces associated with the characteristic foliation $\mathfrak{a}(\mathcal{L})$ of a Hermitian Lie algebroid $\mathfrak{a}\colon \mathcal{L} \rightarrow \mathcal{T}_X$, with the induced inner product (see Theorem~\ref{metric on invariant subspaces}). As an example (following \cite{BP}), we discuss a standard Hermitian Lie algebroid structure on the cotangent sheaf $\mathcal{L} := \Omega^1_X$ of $X$, where $X$ is the $3$-dimensional complex Euclidean space $\mathbb{C}^3$ equipped with a canonical Poisson manifold structure. For a certain $\mathcal{L}$-\emph{invariant subspace} $Y := V(\mathcal{I})$, i.e., $\mathfrak{a}(\mathcal{L})(\mathcal{I}) \subset \mathcal{I}$, we obtain a canonically induced Hermitian Lie algebroid structure on its tangent sheaf $\mathcal{T}_Y$ (see Example~\ref{Foliation with inner product}).

In \cite{JLH} and \cite{GS}, the de Rham cohomologies of (smooth) leaves and leaf spaces associated with smooth foliations on a Riemannian manifold were studied. Here, we present an analogous study in the context of complex geometry, focusing on singular foliations on Hermitian manifolds and their leaves or orbits. To develop a cohomological aspects, we employ associated hypercohomologies, which coincide with the canonical singular cohomology in non-singular cases via the analytic de Rham theorem (see Proposition \ref{cohomology of a leaf}). Building on this cohomology framework, we introduce equivariant analytic de Rham cohomology and extend it to define \emph{equivariant logarithmic de Rham cohomology} (see Section \ref{equiv log de Rham}). To achieve this, we first define $G$-analytic subspaces of a complex manifold with a complex Lie group $G$-action (see Section \ref{G-analytic spaces}). We then show $G$-invariant property and consider the hypercohomology of the associated orbit space, following the theory of Lie-Rinehart cohomology for quotients of singularities by finite groups as developed in \cite{EE}.

In Section \ref{Riemannian geometry for some $a$-spaces}, we consider complex Riemannian geometry of complex affine varieties and Hermitian metrics on analytic spaces. Then in Section  \ref{Hermitian Lie algebroids}, we define Hermitian metrics on holomorphic Lie algebroids and study the associated characteristic foliations with induced metrics. In Section \ref{The analytic de Rham cohomology of spaces associated with Characteristic Foliation}, we consider analytic de Rham cohomologies for leaves, leaf space of a holomorphic foliation on a Hermitian manifold. In Section \ref{Equivarient Lie algebroid cohomology}, we consider equivariant Lie algebroid cohomology in the holomorphic hermitian metric set up.

\section{Riemannian geometry over topological ringed spaces} \label{Riemannian geometry for some $a$-spaces} Riemannian geometry on smooth manifolds is a classical topic in differential geometry, where a Riemannian metric on a smooth manifold $M$ arises from an inner product on its tangent bundle (see \cite{RS, SM, CW}).
 More generally, Lie algebroids over smooth manifolds equipped with an inner product were studied in \cite{MB}, where such structures are referred to as Riemannian metrics on Lie algebroids. In \cite{VP}, Riemannian manifolds were investigated from the viewpoint of Lie--Rinehart algebras, generalizing the notion of an inner product on the $C^\infty(M)$-module $\mathfrak{X}(M)$ of vector fields on $M$ to general $R$-modules over a $\mathbb{K}$-algebra $R$. These $R$-modules are naturally related to the Lie--Rinehart algebra $\operatorname{Der}_{\mathbb{K}}(R)$ of derivations of the $\mathbb{K}$-algebra $R$ (see \cite{GR,JH}).

Here we consider an analogue of these notions, in the context of complex (analytic or algebraic) geometry. For that we consider some particular kind of locally ringed spaces, namely analytic spaces and algebraic varieties etc., considered as topological ringed spaces \cite{JV} or algebraic spaces \cite{AA, Ab}.
\subsection{Complex Riemannian Geometry of complex affine varieties}  \label{Geometry of affine variety}
The notion complex Riemannian geometry appears as an analogue of the Riemannian geometry of smooth manifolds in the holomorphic context \cite{Le}.
  Riemannian Rinehart spaces was introduced in \cite{VP,VP3}, and some of the associated quotient Rinehart spaces viewed as an algebraic generalization of submanifold theory of Riemannian geometry. 
 
 \subsubsection{Riemannian Rinehart spaces}
 We recall the notion (quotient) Riemannian Rinehart spaces introduced by   Pessers and Veken \cite{VP}  and consider some of the key objects  and describe relationships among them. 
\begin{Def} 
By a Rinehart space, we mean a dual pair $\mathcal{A}$ whose primary and secondary $\mathcal{O}$-module denoted by $\mathfrak{X}$ and $\Omega$ respectively, and for which there exists a derivation  $d \in Der_{\mathbb{K}}(\mathcal{O}, \Omega)$, such that the associated map $\tilde{d} \in Hom_{\mathcal{O}}(\mathfrak{X}, Der_{\mathbb{K}}(\mathcal{O}))$ defined by $\tilde{d}D(f)=df(D)$ for $D \in \mathfrak{X}, ~f \in \mathcal{O}$, turns $\mathfrak{X}$ into a $(\mathbb{K}, \mathcal{O})$-Lie algebra (see \cite{GR, MM}) or equivalently say, $(\mathcal{O}, \mathfrak{X})$ forms a Lie-Rinehart algebra (see \cite{JH, MK}).

Furthermore, the triple $(\mathcal{O}, \mathfrak{X}, \Omega)$ is said to be a Rinehart space if the following conditions hold:
\begin{itemize}
	\item The $\mathcal{O}$-modules $\mathfrak{X}$ and $\Omega$ are finitely generated $($or finite rank$)$,
	\item The $\mathcal{O}$-module $\Omega$ is spanned by the image of the map $d: \mathcal{O} \rightarrow \Omega$,
	\item The pairing between $\mathfrak{X}$ and $\Omega$ is non-degenerate.
\end{itemize}

By a Riemannian metric on a Rinehart space $\mathcal{A}:=(\mathcal{O}, \mathfrak{X}, \Omega)$, we simply mean a symmetric, non-degenerate linear map $\langle \cdot, \cdot \rangle : \mathfrak{X} \otimes_{\mathcal{O}} \mathfrak{X} \rightarrow \mathcal{O}$. We call the pair $(\mathcal{A}, \langle \cdot, \cdot \rangle)$ a Riemannian Rinehart space.

Suppose $\mathcal{A}:=(\mathcal{O}, \mathfrak{X}, \Omega)$ is a Rinehart space and $\mathcal{B}:=(\mathcal{I}, \mathfrak{X}', \Omega')$ is a ideal subpair of $\mathcal{A}$, i.e. $\mathcal{I} \subset \mathcal{O}$ is an ideal, $\mathfrak{X}'$ is a Lie ideal of $\mathfrak{X}$ and $\Omega'$ is an $\mathcal{O}$-submodule of $\Omega$ satisfying $d({\mathcal{I}}) \subset \Omega'$.	Then 
	the quotient pair $\mathcal{A}/{\mathcal{B}}:=(\mathcal{O/I},~ \mathfrak{X/X'},~ \Omega/{\Omega'})$ induces a canonical Rinehart space structure,  known as the quotient Rinehart space of $\mathcal{A}$ by $\mathcal{B}$.
	Thus, the quotient map $\rho : \mathcal{A} \rightarrow \mathcal{A}/{\mathcal{B}}$ becomes a surjective morphism of Rinehart spaces.
\end{Def}

In the following, we consider some of the key objects in the context of Riemannian Rinehart space and describe relationships among them.

Let $\mathcal{O}$ be a $\mathbb{C}$-algebra and $\mathcal{I}$ be a finitely presented ideal (prime ideal) of $\mathcal{O}$.
Consider the standard Rinehart space  $$\mathcal{A}:=(\mathcal{O},~ \mathfrak{X}:=Der_{\mathbb{C}}(\mathcal{O}), ~\Omega:=\Omega_{\mathcal{O}/\mathbb{C}})$$ where $Der_{\mathbb{C}}(\mathcal{O})$ is the space of $\mathbb{C}$-linear derivations over the algebra $\mathcal{O}$ and $\Omega_{\mathcal{O}/\mathbb{C}}$ is the space of K$\ddot{a}$hler differentials of the algebra $\mathcal{O}$ over $\mathbb{C}$ (see \cite{JH}).

 The induced quotient Rinehart space of $\mathcal{A}$ associated with ideal $\mathcal{I}$ is $$\mathcal{A}_{\mathcal{I}}:=(\mathcal{O/I},~ \mathfrak{X}_{\mathcal{I}}:=Der_{\mathbb{C}}(\mathcal{O/I}),~ \Omega_{\mathcal{I}}:= \Omega/{<\{dg|~ g \in \mathcal{I}\}>}).$$ 
For the case of polynomial algebra $\mathcal{O}=\mathbb{C}[x_1,\dots,x_n]$ for some $n\in \mathbb{N}$, the spaces $\mathcal{A}$ and $\mathcal{A}_{\mathcal{I}}$ are associated with geometry of the affine algebraic set (affine variety) $Y:=V(\mathcal{I})$ of zero locus of $\mathcal{I}$ inside the affine space $X :=\mathbb{A}^n$, otherwise we can think it as geometry of affine scheme $Y=Spec(\mathcal{O/I})$ associated with $X=Spec(\mathcal{O})$. We consider the canonical action of K$\ddot{a}$hler differentials on derivations as a dual pair $$\langle \cdot, \cdot \rangle ^* : \Omega \otimes_{\mathcal{O}} \mathfrak{X} \rightarrow \mathcal{O}$$ defined as  $(\omega, D) \mapsto \langle \omega, D \rangle ^* := w(D)$ for any $\omega \in \Omega,~ D \in \mathfrak{X}$. If there is a Riemannian metric $\langle \cdot, \cdot \rangle _X$ on the Rinehart space $\mathcal{A}$,
$$\langle \cdot, \cdot \rangle _X : \mathfrak{X} \otimes_{\mathcal{O}} \mathfrak{X} \rightarrow \mathcal{O}$$ given by $(D_1, D_2) \mapsto \langle D_1, D_2 \rangle _X$ for any $D_1, D_2 \in \mathfrak{X}$, then the triple $(\mathcal{O}, \mathfrak{X}, \Omega)$ forms a Riemannian Rinehart space.
In particular for the case of polynomial algebra, these operations are defined by 
\begin{align}\label{standerd inner product}
	\langle df, \sum\limits_{i=1}^{n}f_i \partial_{x_i} \rangle ^*=\sum\limits_{i=1}^{n}f_i ~ \frac{\partial f}{\partial{x_i}}, ~\text{and} ~ \langle \sum\limits_{i=1}^{n}f_i \partial_{x_i} , \sum\limits_{i=1}^{n}g_i \partial_{x_i}\rangle _X= \sum\limits_{i=1}^{n}f_i ~{g_i},
\end{align}
 where $\partial_{x_i}:=\frac{\partial}{\partial{x_i}}$ and $f, f_i, g_i \in \mathcal{O}$ $(i=1, \dots, n)$ and $df(D)=D(f)$ for $D \in \mathfrak{X}$.
In \cite{VP, VP3}, similar ideas are present in a more algebraic approach.
\subsubsection{Particular Rinehart spaces}
Now we construct some of the important (Riemannian) Rinehart spaces for that first we need to consider certain $\mathcal{O}$-modules given as follows.

Consider the $\mathcal{O}$-module ${\mathfrak{X}}^T:= \{D \in Der_{\mathbb{C}}(\mathcal{O}) \mid D(\mathcal{I})\subset \mathcal{I}\}$, consists with such derivations which induces derivations on $\mathcal{O/I}$. In the case of polynomial algebra, it represents space of (algebraic) vector fields on $X:=\mathbb{A}^n$ that are tangent to the zero locus $Y:=V(\mathcal{I})$.

Consider the $\mathcal{O}$-module $\mathfrak{X}^{\perp}:=\{D \in Der_{\mathbb{C}}(\mathcal{O}) \mid \langle D, \mathfrak{X}^T \rangle \subset \mathcal{I}\}$, consists with derivations whose elements represent vector fields normal to $Y:=V(\mathcal{I})$ in $X:=\mathbb{A}^n$ in the case of polynomial algebra.

Consider the $\mathcal{O}$-module $\Omega^{\perp}:= \{w \in \Omega \mid w= f~dg$ for some $ f \in {\mathcal{O}}, g \in {\mathcal{I}}\}$ consists with some special differential $1$-forms.
Thus, we can express ${\mathfrak{X}}^T$ and $\Omega^{\perp}$ as
${\mathfrak{X}}^T = \{D \in Der_{\mathbb{C}}(\mathcal{O}) \mid \langle \Omega^{\perp}, D \rangle \subset \mathcal{I}\} $ and $ \Omega^{\perp} = Hom_{\mathcal{O}}(\mathfrak{X}^T, \mathcal{I}) \cap \Omega$.

Consider the $\mathcal{O}$-module $\mathfrak{X}^0 := \{D \in Der_{\mathbb{C}}(\mathcal{O})\mid \langle \Omega, D \rangle \subset \mathcal{I}\}$,  consists with derivations such that in the case of polynomial algebra $\mathcal{O}=\mathbb{C}[x_1,\dots,x_n]$ each of this derivations on $\mathcal{O}$ (or vector fields over $X=\mathbb{A}^n$) induces a zero derivation on $\mathcal{O}/{\mathcal{I}}$ (or vector field vanishes on $Y$). 

Consider the $\mathcal{O}$-module $\Omega^0 := \{w \in \Omega \mid \langle w, \mathfrak{X} \rangle \subset \mathcal{I}\}$, it represents the K$\ddot{a}$hler differentials that vanishes on $Y=V(\mathcal{I})$ for the polynomial algebra $\mathcal{O}$. 


Consider the quotient map $\rho:\mathcal{O} \rightarrow \mathcal{O/I}$, a surjective $\mathbb{C}$-algebra homomorphism, provide the homomorphism of Lie-Rinehart algebras (see \cite{JH}) 
$$ \rho_* : \mathfrak{X}^T \rightarrow \mathfrak{X}_{\mathcal{I}}$$ defined by $D \mapsto \tilde{D}$. It is well defined because if we take $[{f}]=[{g}]$ that is $f-g \in \mathcal{I}$ for some $f,g \in \mathcal{O}$ then $D(f-g)\in \mathcal{I}$ for any $D \in \mathfrak{X}^T$ which implies $\tilde{D}([{f}])= [{D(f)}]= [{D(g)}]=\tilde{D}([{g}])$. Note that, if $D \in \mathfrak{X}^0, D' \in \mathfrak{X}^T, f \in \mathcal{O}$ then $[D, D'](f) = D(D'f) - D'(Df)\in \mathcal{I}$. Thus, we have  the Lie ideal $$Ker~\rho_* = \{D \in Der_{\mathbb{C}}(\mathcal{O})~|~ D(\mathcal{O})\subset \mathcal{I}\}=\mathfrak{X}^0 \subset \mathfrak{X}^T.$$ 

Note that the differential $d: \mathcal{O} \rightarrow \Omega$ provides the Lie-Rinehart homomorphism 
\begin{center}
$\tilde{d}: \mathfrak{X}^T \hookrightarrow \mathfrak{X}$ defined as $\tilde{d}D(f)= df(D)$
\end{center}
 for any $f\in \mathcal{O}$ and $D \in \mathfrak{X}^T$. Thus $\mathcal{A}^T:=$$(\mathcal{O}, \mathfrak{X}^T, \Omega)$ is a Rinehart subspace of $\mathcal{A}$. 

The differential $d : \mathcal{O} \rightarrow \Omega$ induces $\bar{d} := d|_{\mathcal{I}} : \mathcal{I}\rightarrow \Omega^{\perp}$, which provides a natural $(\mathbb{C}, \mathcal{I})$-Lie-Rinehart algebra homomorphism  $ \bar{d} : \mathfrak{X}^T \rightarrow Der_{\mathbb{C}}(\mathcal{I})$ (if $f \in \mathcal{O}, g \in \mathcal{I}$ and $ D\in \mathfrak{X}^T$ then $\langle f~ dg, D \rangle = \langle dg, f~D \rangle = f$ $ D(g) \in \mathcal{I})$. Thus we get the Rinehart space $(\mathcal{I}, \mathfrak{X}^T, \Omega^{\perp})$.

\begin{Rem}
The above described Rinehart spaces can be viewed as Riemannian Rinehart spaces by considering the induced structure coming from a Riemannian Rinehart space structure on $(\mathcal{O}, \mathfrak{X}, \Omega)$.
\end{Rem}

Most of the notions discussed so far are based on the framework established in \cite{VP}. These concepts serve as a local model within the realm of algebraic geometry. By globalizing these ideas through sheaf theory, we obtain an analogue of Riemannian geometry for algebraic varieties and, more generally, for schemes. Here, we explore the corresponding global version within the context of complex geometry.

\subsubsection{Holomorpic Riemannian manifolds}
In \cite{DZ,Le,VP2}, the notion of a holomorphic Riemannian manifold, defined as a complex manifold $X$ endowed with a holomorphic Riemannian metric, is considered. A holomorphic Riemannian metric on a complex manifold $X$ is a holomorphic section of the tensor bundle $\otimes^2 T^*X$ (i.e., a holomorphic covariant $2$-tensor) which is symmetric and non-degenerate. In other words, it is a holomorphic field of non degenerate complex quardetic forms on the holomorphic tangent bundle $TX$. 
    
    For example, consider $\mathbb{C}^n$ with the global flat holomorphic Riemannian metric $\sum_{j=1}^{n} (dx_j)^{\otimes 2}$ $($analogue of the inner product (\ref{standerd inner product})$)$ where $\{x_j\}^n_{j=1}$ is the standard coordinate system of $\mathbb{C}^n$. Take any quotient of $\mathbb{C}^n$ by a lattice $(n$-dimensional complex torus is one such space$)$ with induced metric.

    Moving forward, we apply these ideas to develop the foundations of complex Riemannian geometry and Hermitian geometry for analytic spaces. This approach enables the generalization of classical differential geometry for submanifold theory of Riemannian manifolds to the holomorphic setting.

\subsection{Hermitian metrics on analytic spaces} \label{Hermitian metric on analytic spaces }  By following the concepts of Hermitian metric on a complex vector bundle \cite{GH} and Riemannian manifolds as Lie-Rinehart algebras \cite{VP,VP3}, we consider an analogue of Hermitian metric in the holomorphic context for analytic spaces. 

First, we recall some important notions below, which will be useful for the remainder of the article.

\subsubsection{Holomorphic tangent bundle} Utilizing the underlying smooth structure of a complex manifold $X$, on each open set $U\subset X$, consider the $\mathbb{C}$-algebra of complex valued smooth (or $C^\infty$) functions on $U$:
$$C_X^\infty(U):=\{f:U \rightarrow \mathbb{C}~|~ f~\text{is a smooth map}\}.$$
 The induced sheaf of $\mathbb{C}$-valued $C^\infty$ functions on $X$ is denoted by $C_X^\infty$, combining the algebraic structure of $\mathbb{C}$  with the smooth structure of $\mathbb{R}^2$.
The sheaf of derivations over the sheaf of holomorphic functions $\mathcal{O}_X$ is $$\mathcal{D}er_{\mathbb{C}_X}(\mathcal{O}_X):=\{D: \mathcal{O}_X \rightarrow \mathcal{O}_X~|~ D~\text{is a}~\mathbb{C}\text{-linear sheaf homomorphism,}~\text{satisfying the Leibniz rule}\},$$ isomorphic to the sheaf of sections $\Gamma_X(TX)=:\mathfrak{X}_X$ of the holomorphic tangent bundle $TX$ over $X$ (i.e., the sheaf of holomorphic vector fields on $X$). It has a canonical $\mathbb{C}_X$-Lie algebra and left $\mathcal{O}_X$-module structure, known as the tangent sheaf $\mathcal{T}_X:=(\mathcal{D}er_{\mathbb{C}_X}(\mathcal{O}_X), [\cdot, \cdot]_c)$ of $X$.
The complex conjugate of the $\mathbb{C}$-vector space $T_xX$ (tangent space at $x \in X$) is denoted by $\overline{T_xX}$, and the associated conjugate bundle of the holomorphic bundle $TX$ is denoted by $\overline{TX}$. We consider the smooth $\mathbb{C}$-vector bundle $ (TX \otimes \overline{TX})^*$, whose fiber is $$(T_xX \otimes_{\mathbb{C}} \overline{T_xX})^*:= Hom_{\mathbb{C}} (T_xX \otimes_{\mathbb{C}} \overline{T_xX}, ~\mathbb{C})~\text{at}~x \in X.$$ 

\subsubsection{Hermitian metric on a complex manifold $(X, \mathcal{O}_X)$} Let $\{(U_i, \phi_i)\}_{i \in \mathbb{N}}$ be a holomorphic atlas of a complex manifold $X$ (dimension of $X$ is $n$), i.e.  $\phi_i: U_i \rightarrow V$ is a biholomorphism with some proper open subset $V$ of $\mathbb{C}^n$. Then we can induce an inner product locally on the holomorphic tangent bundle $TX$ from the standard inner product on $\mathbb{C}^n$, and globalize it by using a $C^\infty$-partition of unity associated to the open cover of $X$. This provides us a $\mathbb{C}$-sesquilinear, conjugate-symmetric, positive definite map
$$h_x: T_xX \times T_xX \rightarrow \mathbb{C}$$
on each holomorphic tangent space $T_xX$ (at $x \in X)$, varies smoothly on $X$. 

It can be viewed as smoothly varying positive definite Hermitian form $($or Hermitian inner product$)$ 
$$\langle \cdot, \cdot \rangle _x : T_xX \otimes_{\mathbb{C}} \overline{T_xX} \rightarrow \mathbb{C}~,$$ where $\langle v, c~ w \rangle _x:= h_x(v, \bar{c}w) =c~ h_x(v,w) =c \langle v, w \rangle _x$, for $c \in \mathbb{C}$,~ $v,w \in T_xX$ (i.e. $\langle \cdot, \cdot \rangle _x \in (T_xX \otimes_{\mathbb{C}} \overline{T_xX})^*)$. It induces the $C^\infty$-map $($which forms a smooth global section of the vector bundle $(TX \otimes \overline{TX})^*$ on $X)$
\begin{center}
	$\langle \cdot, \cdot \rangle : X \rightarrow (TX \otimes \overline{TX})^*$ given by\\
	\hspace{-1.5 cm}	$x \mapsto \langle \cdot, \cdot \rangle _x$
\end{center}
for any $x \in X$. Such a smoothly varying positive-definite Hermitian form on each fibre $T_xX$ of the holomorphic tangent bundle $TX$ is known as a Hermitian metric on $X$. 
Thus, on an open set $U_i$ associated with a chart, we get a map given by the positive definite Hermitian form
\begin{align}\label{local inner product}
	\langle \cdot, \cdot \rangle _{U_i} : \mathcal{T}_X(U_i) \times  {\mathcal{T}_X(U_i)} \rightarrow C_X^\infty(U_i),
\end{align}
i.e. it is an $\mathcal{O}_X(U_i)$-sesquilinear, conjugate-symmetric and positive definite map.

Now, by using a $C^\infty$-partition of unity $\{f_i\}_{i \in \mathbb{N}}$ for the open cover $\{U_i\}_{i \in \mathbb{N}}$, we can form a Hermitian inner product on $X$ as $$\langle \cdot, \cdot \rangle_X := \sum_{i \in \mathbb{N}} f_i \langle \cdot, \cdot \rangle_{U_i}$$ i.e., in the global situation by considering respective sheaves, we get the $\mathcal{O}_X$-sesquilinear map 
$$\langle \cdot, \cdot \rangle_X: \mathcal{T}_X \times {\mathcal{T}_X} \rightarrow C^\infty_X$$ given by the local Hermitian forms as
\begin{align} \label{global inner product}
   \langle gD, g'D' \rangle _X= \sum_{i \in \mathbb{N}} f_i~g \bar{g}' \langle D|_{U_i}, D'|_{U_i}\rangle _{U_i}, 
\end{align}
for $D,D' \in \mathcal{T}_X \cong \mathfrak{X}_X$ and $g, g' \in \mathcal{O}_X$. 
Additionally, the map satisfies the properties of being conjugate-symmetric and positive definite on the space of sections. Note that, here the sheaf $\mathcal{T}_X \times {\mathcal{T}_X}$ is the direct product $($or direct sum$)$ of $\mathcal{O}_X$-modules.

A Hermitian manifold is  a complex manifold $X$ with
a Hermitian metric on $X$, which is a positive definite hermitian form on the sheaf of holomorphic vector fields $\mathfrak{X}_X$ of $X$.
\begin{Rem}\label{standard inner product}
	On a chart $(U,(x_1, \dots, x_n))$ of a complex manifold $(X, \mathcal{O}_X)$, a Hermitian metric is given by 
    $$\sum_{1 \leq i,j \leq n} h_{ij} dx_i \otimes d\bar{x_j},$$ where each function $h_{ij}$ is determined by the Hermitian inner product on $T_xX$ as $h_{ij}(x)=\langle \partial_{x_i}, \partial_{x_j} \rangle_x$ for $x \in U$. The standard Hermitian metric on an open subset $U$ of $\mathbb{C}^n$ is given by 
    $$\langle \sum^n_{i=1} f_i \partial_{x_i}, \sum_{j=1}^{n} g_j \partial_{x_j} \rangle_U= \sum^n_{i=1}f_i  \bar{g}_i,~\text{for all}~f_i, g_i \in \mathcal{O}_{\mathbb{C}^n}(U),~ \text{where}~i=1, \dots, n.$$
\end{Rem}

	The concept of a Hermitian metric on a complex manifold, i.e. a smoothly varying Hermitian inner product on each $($holomorphic$)$ tangent space, extends to the idea of considering a Hermitian metric on a holomorphic vector bundle, considering smoothly varying fiberwise Hermitian inner product (see \cite{GH}).

		We extend the idea of complex conjugate of a holomorphic vector bundle on a complex manifold $X$ (see \cite{GH}) by considering complex conjugate on its sheaf of sections of an $\mathcal{O}_X$-module.

\subsubsection{Lie algebroids over analytic spaces}
We recall the definition of Lie algebroids in the algebro-geometric language in order to unify various geometric notions required for this article (see \cite{MK, BP, BRT, AA, PT, TS}). 
\begin{Def}  
A Lie algebroid $(\mathcal{L}, [\cdot, \cdot], \mathfrak{a})$ over an analytic space $(X, \mathcal{O}_X)$ consists of a coherent $\mathcal{O}_X$-module $\mathcal{L}$ with a $\mathbb{C}_X$-Lie algebra structure on it given by $[\cdot, \cdot]$, and a homomorphism $$\mathfrak{a} : (\mathcal{L}, [\cdot, \cdot]) \to (\mathcal{D}er_{\mathbb{C}_X}(\mathcal{O}_X), [\cdot, \cdot]_c)$$ of $\mathcal{O}_X$-modules and $\mathbb{C}_X$-Lie algebras, called the anchor map. The map $\mathfrak{a}$ satisfies the Leibniz rule:
$$[D, f D'] = f [D, D'] + \mathfrak{a}(D)(f) D', ~\text{for}~ f \in \mathcal{O}_X ~\text{and}~D, D' \in \mathcal{L}.$$  
That is, $(\mathcal{L}, [\cdot, \cdot], \mathfrak{a})$ is a coherent sheaf of $(\mathbb{C}_X, \mathcal{O}_X)$-Lie-Rinehart algebras, which we denote simply by $\mathcal{L}$.

Let $(\mathcal{L},[\cdot,\cdot],\mathfrak{a})$ and $(\mathcal{L}',[\cdot,\cdot]',\mathfrak{a}')$ be Lie algebroids over analytic spaces $(X, \mathcal{O}_X)$ and over $(X', \mathcal{O}_{X'})$ respectively, such that $\phi: \mathcal{O}_X \rightarrow \mathcal{O}_{X'}$ is a $\mathbb{C}_X$-algebra homomorphism.
			A homomorphism of the Lie algebroids $$ \psi:  (\mathcal{L},[\cdot,\cdot],\mathfrak{a}) \rightarrow (\mathcal{L}',[\cdot,\cdot]',\mathfrak{a}')$$  is a sheaf homomorphism of Lie-Rinehart algebras (see \cite{JH} for the local description), i.e. 
			 \begin{itemize}
				\item $\psi: \mathcal{L} \rightarrow \mathcal{L}'$ is a morphism of $\mathcal{O}_X$-modules where $\mathcal{O}_X$ acts on $\mathcal{L}'$ via $\phi$,
				\item $\psi :  (\mathcal{L},[\cdot,\cdot])\rightarrow (\mathcal{L}',[\cdot,\cdot]')$ is a $\mathbb{C}_X$-Lie algebra homomorphism,
				\item compatibility:  $\mathfrak{a}' (\psi (D)) (\phi(f))= \phi(\mathfrak{a}(D)(f))$, for $f \in \mathcal{O}_X$ and $D \in \mathcal{L}$.
			\end{itemize}

\end{Def} 
\begin{Exm} For a complex manifold $(X,\mathcal{O}_X)$, consider the analytic subspace  $(Y:=V(\mathcal{I}),~\mathcal{O}_Y:=\mathcal{O}_X/{\mathcal{I}})$ defined by an ideal sheaf $\mathcal{I}\subset \mathcal{O}_X$.
The sheaf of logarithmic derivations  of $X$ along $Y$ is given by
$$\mathcal{T}_X(-log~Y):=\{D \in \mathcal{T}_X~|~ D(\mathcal{I}) \subset \mathcal{I}\},$$  has a canonical Lie algebroid structure. Also, the tangent sheaves $\mathcal{T}_X:=\mathcal{D}er_{\mathbb{C}_X}(\mathcal{O}_X)$,~  $\mathcal{T}_Y:=\mathcal{D}er_{\mathbb{C}_Y}(\mathcal{O}_Y)$ of $X$ and $Y$ respectively, have the standard Lie algebroid structure (see \cite{BP,AA}).
The quotient map 
\begin{align} \label{rho}
  \rho: \mathcal{O}_X\rightarrow \mathcal{O}_X/{\mathcal{I}}=\mathcal{O}_Y~~\text{is given by}~~f\mapsto [f],  
\end{align}
 a surjective $\mathbb{C}_X$-algebra homomorphism. Consider the canonical surjective Lie algebroid morphism $$\tilde{\rho}: \mathcal{T}_X(-log ~Y)\rightarrow \mathcal{T}_Y~~\text{defined by}~~\tilde{\rho}(g D)=[g] \tilde{D},$$ where $\tilde{D}([f])=[D(f)]$, for $f, g \in \mathcal{O}_X$ and $D \in \mathcal{T}_X(-log~ Y)$. 
Thus, $\mathcal{K}er$ $\tilde{\rho}=\mathcal{I}\mathcal{T}_X=:\mathcal{T}^0_X$.
\end{Exm}

\begin{Rem}	
	A locally free Lie algebroid over a complex manifold is equivalent to a holomorphic Lie algebroid \cite{BP}, where a holomorphic Lie algebroid  is defined as an analogue of smooth Lie algebroids \cite{KM} in the holomorphic category. These are the classical cases, where we have certain Lie-Rinehart algebra structures on the space of sections of a smooth or holomorphic vector bundle \cite{AA}.
	\end{Rem}
We consider Hermitian metrics on Lie algebroids over analytic spaces, as a generalization of Hermitian manifold. For that first we prove some results related to Hermitian manifold. Then we define the analogue for Lie algebroids over analytic spaces with Hermitian metrics and consider some of its special cases.
\begin{Thm}
	Let $(X, \langle \cdot, \cdot \rangle)$ be a Hermitian manifold and $U$ an open set contained in some chart of $X$. Let $D \in \mathcal{T}_X(U)$ such that $\langle D, D' \rangle \in \mathcal{I}(U)$ for all $D' \in \mathcal{T}_X(U)$, where $\mathcal{I}\subset \mathcal{O}_X$ is an ideal sheaf. Then $D \in \mathcal{T}^0_X(U)$, and moreover for any global section $\tilde{D}$ of $\mathcal{T}_X$ if locally this result holds, then it implies that $\tilde{D} \in \mathcal{T}^0_X (X)$. 
\end{Thm} 

\begin{proof}
	Let $ x_1, \dots, x_n $ be the co-ordinate functions on the chart $U$. Then $D \in \mathcal{T}_X(U)$ is expressed as
    $$D = \sum^n_{i=1}g_i \partial_{x_i}~\text{for some}~g_1, \dots, g_n \in \mathcal{O}_X(U).$$ 
    Take $D' = \partial_{x_j}$ and thus we get $g_j = \langle \sum^n_{i=1}g_i \partial_{x_i}, \partial_{x_j} \rangle \in \mathcal{I}(U)$, for $j = 1, \dots, n$. 
    Notice that,  $$~~\text{if}~~D(\mathcal{O}_X(U)) \subset \mathcal{I}(U),~~\text{then}~~g_j = D(x_j)\in \mathcal{I}(U)$$ and  $$~~\text{if}~~D \in \mathcal{I}(U)\mathcal{T}_X(U),~~\text{then}~~D(f) = \sum^n_{i=1} g_i  \partial_{x_i}(f)\in \mathcal{I}(U)~~\text{for all}~~f \in \mathcal{O}_X(U).$$  Hence, $D \in \mathcal{I}(U)\mathcal{T}_X(U) = \mathcal{T}^0_X(U)$.
	If $\tilde{D} \in \mathcal{T}_X(X)$ satisfy the above condition locally, Then, for each point $x \in X$ there exists an open set $U_x$ contained in some chart such that $\tilde{D} = D$ on $U_x \subset U$. Thus, the stalk of $\tilde{D}$ at $x$ is $D_x$, and hence $D_x \in \mathcal{T}^0_{X,x}$. Therefore, $\tilde{D} \in \mathcal{T}^0_X(X)$.
\end{proof}

\subsubsection{The induced metric on analytic subspaces or divisors}
If $ \mathcal{I} \subset \mathcal{O}_X$ is a principal ideal sheaf (implies $Y:=V(\mathcal{I})$ is a hypersurface singularities), then for each point $x \in X$ there is an open neighborhood  $U$ (contained in some chart) around $x$ such that $\mathcal{I}(U)$ is generated by some $g \in \mathcal{O}_X(U)$, i.e. $\mathcal{I}(U)=\langle g \rangle$. Then for any $D \in \mathcal{T}_X(U)$, define the derivation on $\mathcal{O}_X(U)$ by
\begin{align}\label{along tangent} 
	D^T : = D - \langle D, \nabla g \rangle \nabla g
\end{align}
where $$\nabla g = \sum^n_{j=1} \partial_{x_j}(g)~\partial_{x_j} \in \mathcal{N}_{Y/X}(U)$$ is the gradient of $g \in \mathcal{O}_X(U)$ and $\mathcal{N}_{Y/X}$ is the sheaf of normals on $Y$ inside $X$. To show $D^T=D-\langle D, \nabla g \rangle \nabla g \in \mathcal{T}_X(-log~Y)(U)$, we need to check $D^T(f ~ g) \in \mathcal{I}(U)$ for any $f \in \mathcal{O}_X(U)$. Notice that we have the following identity 
\begin{align*}
	D^T(f ~ g)
	=&~ D(f ~ g)- D(g) ~ \nabla(g)(f ~ g)\\
	=&~ D(f) ~ g + f ~ D(g) - D(g) ~ (\nabla g(f) ~ g + f ~ \nabla(g)(g)) \\
	=&~(D- D(g) \nabla g)(f) ~ g + D(g) ~ f ~ (1- \langle \nabla g, \nabla g \rangle).
\end{align*}
In the above identity the first term is in the ideal $\langle g \rangle = \mathcal{I}(U)$ and the second term also is in $\mathcal{I}(U)$ if the following condition holds:
\begin{align}\label{cond 1} 
	1- \langle \nabla g, \nabla g \rangle \in \mathcal{I}(U).
\end{align}
 In particular, when $Y$ is a complex submanifold of $X$, of codimension $1$, we choose an atlas of $X$ such that $Y$ is locally represented as $\{x_i=0\}$ for some coordinate function $x_i$ in some chart of $X$. In that case we have the identity $$ \langle \nabla g, \nabla g \rangle = \nabla g(g)= \sum_{j=1}^{n}  {\partial_{x_j}}(g) \cdot {\partial_{x_j}}(g)=1$$ (putting $g=x_i$) and so the required condition holds. 
But in general, the condition (\ref{cond 1}) does not hold for hypersurface singularities. Like when we consider the normal crossing divisor $Y=V(xy)$ in $X=\mathbb{C}^2$, we face difficulties. In this case the principal ideal sheaf is $\mathcal{I}=\langle xy \rangle$ and $\mathcal{O}_X$ is the sheaf of holomorphic functions. Thus, the sheaf of logarithmic derivations is $\mathcal{T}_X(- log~ Y)=\langle \{x \partial_x, y \partial_y\} \rangle$ and the sheaf of normals $\mathcal{N}_{Y/X}=\langle y \partial_x + x \partial_y \rangle$. Since here $g=xy$, we obtain the following: 
$$ \langle \nabla g, \nabla g \rangle = \nabla g(g)= x^2+y^2,~ \text{and}~ 1- \langle \nabla g, \nabla g \rangle = 1-\nabla g(g)= 1-(x^2+y^2) \notin \langle \{xy\} \rangle$$

Therefore if the condition holds, we get the Lie-Rinehart algebra homomorphisms (see (\ref{along tangent}))
\begin{align} \label{TX to TY}
 \mathcal{T}_X(U) \overset{(\cdot)^T _U}{\rightarrow} \mathcal{T}_X(-log~ Y)(U)\overset{\rho_{*U}}{\rightarrow} \mathcal{T}_Y(Y \cap U)   
\end{align}
\begin{center}
\hspace{-1 cm} $D \longmapsto D^T \longmapsto \widetilde{D^T}:=\bar{D}$,
\end{center}
where $\rho_*$ is the induced map from the map $\rho$ (\ref{rho}).
Hence, in each of the charts we have these morphisms which are compatible with any open subsets of a chart. Induced morphisms provides stalkwise morphisms of Lie-Rinehart algebras, sheafifying it we get the homomorphisms of Lie algebroids 
\begin{align} \label{pi}
   \mathcal{T}_X \overset{(\cdot)^T }{\rightarrow} \mathcal{T}_X(-log~ Y) \overset{\rho_{*}}{\rightarrow} \mathcal{T}_Y, 
\end{align}
locally defined by the above way. We denote the morphism $\rho_* \circ (\cdot)^T$ by $\pi$.

\smallskip
\begin{Nota} \label{Notation}
    We use the fact that $\mathcal{I} \subset \mathcal{O}_X \subset C^\infty_X$ holds for a complex manifold $X$. Notice that, the quotient sheaf $C^\infty _Y:= C^\infty _X/{\mathcal{I}}$ on $Y$, forms an $\mathcal{O}_X$-module. Thus, the quotient map 
    \begin{align}\label{rho'}
        \rho': C_X^\infty \rightarrow C_X^\infty/{\mathcal{I}}
    \end{align}
    $($the canonical $\mathcal{O}_X$-module homomorphism$)$ satisfies $\rho'|_{\mathcal{O}_X}=\rho$, useful to define an inner product on $\mathcal{T}_Y$.

From now, we simply write $X$ to denote a Hermitian manifold $(X, \langle \cdot, \cdot \rangle_X)$, and denote the restriction $\langle \cdot, \cdot \rangle_X|_U$ of $\langle \cdot, \cdot \rangle_X$ on $U$ by  $\langle \cdot, \cdot \rangle_U$ for an open subset $U$ of $X$.
\end{Nota} 
\begin{Thm}\label{induced inner product}
	Let $X$ be a Hermitian manifold and $Y$ an analytic subspace of X satisfying the condition (\ref{cond 1}). Then $Y$ inherits from $X$, an $\mathcal{O}_Y$-sesquilinear form $\langle \cdot, \cdot \rangle_Y$ which is locally determined by the equation: for any $\bar{D}, \bar{D}' \in \mathcal{T}_Y(U)$, $\langle \bar{D}, \bar{D'}\rangle _{Y\cap U}:= \rho'(U) (\langle D, D'\rangle _U)$ for some $D, D' \in \mathcal{T}_X(U)$ where $U$ is an open set contained in a chart of $X$, such that $\pi(U)(D) = \rho_*(U)(D^T) = \bar{D}$ and $\pi(U)(D') = \rho_*(U)(D'^T) = \bar{D'}$ holds.
\end{Thm}
\begin{proof}
	Well definedness: Let $D, D'' \in \mathcal{T}_X(U)$ satisfying the identity (see (\ref{pi}))
    $$\pi(U)(D) = \rho_*(U)(D^T) = \bar{D} = \rho_*(U)(D''^T) = \pi(U)(D'').$$ Then $\pi(U)(D-D'') = \bar{0}$, i.e. $D-D'' \in \mathcal{T}^0_X(U)$ and hence $\langle D-D'', D' \rangle _U \in \mathcal{I}(U)$. Thus, 
    $$\rho'_U(\langle D-D'', D' \rangle _U) = \bar{0}~~\text{or}~~\rho'_U(\langle D, D' \rangle _U)= \rho'_U(\langle D'', D' \rangle _U),$$ where $\rho'_U=\rho|_U$ (see (\ref{rho}) and (\ref{rho'})). Therefore, for an open set $U_i$ associated with a chart we have the induced $\mathcal{O}_Y|_{Y \cap U_i}$-bilinear form $\langle \cdot, \cdot \rangle_{Y \cap U_i}$ on $Y \cap U_i$ as defined in the statement. Now using the $C^\infty$-partition of unity $\{f_i|_{Y \cap U_i}\}$ of $Y$, we can form the induced  $\mathcal{O}_Y$-sesquilinear form  $\langle \cdot, \cdot\rangle _Y$ on $Y$ as $$\langle \cdot, \cdot \rangle _Y := \sum_{i} f_i|_{Y \cap U_i} \langle \cdot, \cdot\rangle _{Y \cap U_i},$$ 
	induces an $\mathcal{O}_Y$-sesquilinear sheaf homomorphism 
	$$ \langle \cdot, \cdot\rangle _{Y}: \mathcal{T}_Y \times {\mathcal{T}_Y} \rightarrow C^\infty _Y.$$
\end{proof}
Note that, Theorem \ref{induced inner product} is a global counterpart of a result that appeared in \cite{VP,VP3}.
\begin{Rem}\label{non degenerate property}
The conjugate-symmetry property of the induced $\mathcal{O}_Y$-sesquilinear form $\langle \cdot, \cdot \rangle_{Y}$ on $\mathcal{T}_Y$ follows from the conjugate-symmetry property of the inner product on $\mathcal{T}_X$ as follows.
Notice that, for $\bar{D}, \bar{D}' \in \mathcal{T}_Y$ we have the identity
$$\langle \bar{D}, \bar{D'}\rangle _{Y\cap U}:= \rho'_U (\langle D, D'\rangle _U) = \rho'_U (\overline{\langle D', D\rangle _U}) = \overline{\langle \bar{D'}, \bar{D}\rangle}_{Y\cap U}.$$ 
\label{cond 2} But the positive definite property of $\langle \cdot, \cdot\rangle _{Y}$ does not induces naturally: if for any $\bar{D} \in \mathcal{T}_Y(Y \cap U)$ such that $$\langle \bar{D}, \bar{D'}\rangle _{Y\cap U} = \bar{0}~\text{holds for all}~\bar{D'} \in \mathcal{T}_Y(Y \cap U),~\text{then we obtain}~ \langle D, D'\rangle _U \in \mathcal{I}(U), ~\text{for all}~D' \in \mathcal{T}_X(-log~ Y)(U).$$ This does not implies $D \in \mathcal{T}^0_X(U)$ always, but here  it holds using condition (\ref{cond 1}) and thus $\bar{D} =0$.

	Thus, in this case we get an inner product structure on $\mathcal{T}_Y$, we call it the induced metric on $Y$.
\end{Rem}
\begin{Cor}
	For a Hermitian manifold $X$, by considering the restriction of the Hermitian inner product $\langle \cdot, \cdot \rangle_X$ on $\mathcal{T}_X$ to $\mathcal{T}_X(-log~ Y)$ where $Y$ is an analytic subspace of $X$ satisfying condition (\ref{cond 1}), provides a Hermitian inner product on $\mathcal{T}_X(-log~ Y)$ and $\mathcal{T}_Y$.
\end{Cor}
\subsection{Covariant connections} \label{connection} For a locally free $\mathcal{O}_X$-module $\mathcal{E}$ of rank $r$ over a complex manifold ($X, \mathcal{O}_X$) for some $r \in \mathbb{N}$, with a holomorphic atlas $\{(U_i, \phi_i)\}_{i \in \mathbb{N}}$ of $X$, we have a canonical covariant derivative connection on $\mathcal{E}$, given as follows. Consider an $\mathcal{O}_X(U_i)$-linear map (associated with an open set $U_i$) 
$$\nabla^{i}: \mathcal{T}_X(U_i) \rightarrow End_{\mathbb{C}}(\mathcal{E}(U_i))$$
satisfies the Leibniz rule
	$\nabla^{i}_{D'}(fs)  = f\nabla^{i}_{D'}(s) + D'(f)s$, where $\nabla^{i}_{D'}:=\nabla^{U_i}(D')$,
for sections $f \in \mathcal{O}_X(U_i)$, $D' \in \mathcal{T}_X(U_i)$ and $s \in \mathcal{E}(U_i)$. There is a covariant  connection (see \cite{SM}) satisfying the local descriptions  
\begin{align} \label{covariant derivative locally}
   \nabla^{i}_D( \sum^r_{i=1} g_i s_i) = \sum^r_{i=1} D(g_i) s_i, 
\end{align}
where $D \in \mathcal{T}_X(U_i)$ and $\{s_1, \dots, s_r\}$ is a basis of $\mathcal{E}(U_i)$ (without loss of generality we consider a common open cover $\{U_i\}_i$ of $X$ where both $\mathcal{T}_X(U_i)$ and $\mathcal{E}(U_i)$ are free $\mathcal{O}_X(U_i)$-modules of finite rank).
Now by using a $C^\infty$-partitions of unity $\{f_i\}_i$ on this open cover of $X$, we form the $C^{\infty}_X$-module homomorphism 
$$\nabla: {C^{\infty}_X} \otimes_{\mathcal{O}_X} \mathcal{T}_X \rightarrow {C^{\infty}_X} \otimes_{\mathcal{O}_X} \mathscr{E}nd_{\mathbb{C}_X}(\mathcal{E}),$$  defined by $\nabla := \sum_{i} f_i ~ \nabla^{i}$, i.e. $\nabla^U_D(s)= \sum_{i} f_i ~ \nabla^{i}_D(s)|_U$ where $\nabla^U:=\nabla|_U$,
for $D \in \mathcal{T}_X|_U$ and $s \in \mathcal{E}|_U$.

In particular for $\mathcal{E}=\mathcal{T}_X$, we say that a covariant derivative $\nabla$ is symmetric if we have 
  $$\nabla^U_D D' - \nabla^U_{D'} D = [D,D']_U,$$
for all $D, D' \in \mathcal{T}_X|_U$,
and is compatible with an inner product $\langle \cdot,\cdot \rangle_X$ (see (\ref{local inner product}) and (\ref{global inner product})) if  we have 
$$D (\langle D',D'' \rangle_U) = \langle \nabla^U _D D', D'' \rangle _U + \langle D', \nabla^U _D D'' \rangle _U,$$
for sections $D, D', D'' \in \mathcal{T}_X|_U$.
We say that a covariant derivative connection is a Levi-Civita connection if it is both symmetric and compatible with the inner product on the tangent bundle. The canonical covariant derivative on a complex manifold  satisfies these properties and is thus a Levi-Civita connection.

\begin{Rem}
	If  $F \subset TX$ is an involutive subbundle of the tangent bundle of a Hermitian manifold $X$, then the corresponding normal bundle of $F$ in $TX$ provides a locally free sheaf of $\mathcal{O}_X$-module (formed by its sheaf of sections). Thus it has a canonical covariant connection as described above. Later in Section \ref{The analytic de Rham cohomology of spaces associated with Characteristic Foliation},
     we use this concept to define a version of the Bott connection \cite{JH}.
\end{Rem} 
\begin{Rem} The concept of covariant connections extends to Lie algebroid connections in (smooth real)  differential and complex geometric settings (see \cite{KM, RF,BP}).
	Consider the Lie algebroid $\mathcal{L}=\mathcal{T}_X(- log~ Y)$ for an analytic subspace $Y$ in $X$. If an $\mathcal{L}$-connection exists, it is called a logarithmic connection which is a meromorphic connection with simple poles along the divisor $Y$ \cite{BP, Ab}.
\end{Rem}
\begin{Thm} \label{connection on $Y$}
	Let $X$ be a Hermitian manifold with a Levi-Civita connection $\nabla$ and $Y$ an analytic subspace of $X$ with the induced metric. Then $Y$ inherits a Levi-Civita connection $\bar{\nabla}$ from the connection $\nabla$ on $X$, determined by the equation:
	$${\bar{\nabla}^{U}}_{\bar{D}}(\bar{D'}) = \pi_U(\nabla^U _D D'),$$
 for all $\bar{D} = \pi_U(D),~ \bar{D'} = \pi_U(D')$, where $\pi_U:=\pi|_U$ and $U$ is an open set of $X$ (see (\ref{pi})).
\end{Thm}
\begin{proof}
	Well definedness: Let $\pi_U(D) = \bar{D} = \pi_U(D_1)$. Then $D'':=D-D_1 \in \mathcal{T}^0_X|_U$, hence $\nabla^U_{D''}D' \in \mathcal{T}^0_X|_U$. 
    $$\text{This happens, since}~ D''=\sum^k_{i=1}g_i~ D_i~\text{for some}~g_1, \dots, g_k \in \mathcal{I}|_U,~\text{thus}~\nabla^U_{D''}D' = \sum^k_{i=1}g_i~\nabla^U_{D_i}D'\in \mathcal{I}|_U \mathcal{T}_X|_U.$$ 
    Therefore, $\pi_U(\nabla^U_{D''}D') = \bar{0}$ holds, and since we obtain the identity 
    $${\bar{\nabla}^{U}}_{\bar{D''}}(\bar{D'}) = \pi_U(\nabla^U_{D''}D') = \pi_U(\nabla^U_{D}D' -\nabla^U_{D_1}D'),~\text{thus}~\pi_U(\nabla^U_{D}D') = \pi_U(\nabla^U_{D_1}D').$$ 
    Now, let $\pi_U(D') = \bar{D'} = \pi_U(D'_1)$. Then 
    $D'':=D'-D'_1 \in \mathcal{T}^0_X|_U,$ as of the above form. Without loss of generality we can take $D \in \mathcal{T}_X(-log~ Y)|_U$, because $\pi_U(D)=\pi_U(D^T)$ for every $D \in \mathcal{T}_X|_U$ (see (\ref{TX to TY})) and thus 
    $$\pi_U(\nabla^U_{D}D'') = \bar{0},~\text{follows from the fact}~\nabla^U_{D}D'' = \sum^k_{i=1}D(g_i)D_i  \in \mathcal{I}|_U \mathcal{T}_X|_U=\mathcal{T}^0_X|_U.$$  
   Since we have ${\bar{\nabla}^{U}}_{\bar{D}}(\bar{D''}) = \pi_U(\nabla^U_{D}D'') = \pi_U(\nabla^U_{D}D' -\nabla^U_{D}D'_1),$ hence $\pi_U(\nabla^U_{D}D') = \pi_U(\nabla^U_{D}D'_1)$.
Therefore, the induced connection $\bar{\nabla}$ on $Y$ is well defined. 

\smallskip	
	
\noindent	Symmetry: 
	\begin{align*}
		\hspace{-2 cm} {\bar{\nabla}^{U}}_{\bar{D}}(\bar{D'}) - {\bar{\nabla}^{U}}_{\bar{D'}}(\bar{D})
		=&~ \pi_U(\nabla^U _D D' - \nabla^U _{D'} D) \\
		=&~ \pi_U([D, D']_U) \\
		=&~ [\pi_U(D), \pi_U(D')]_{Y \cap U} \\
		=&~ [\bar{D}, \bar{D'}]_{Y \cap U}.
	\end{align*}
	Compatibility: 
	\begin{align*}
		& \langle {\bar{\nabla}^{U}}_{\bar{D}}(\bar{D'}), \bar{D''} \rangle _U + \langle \bar{D'}, {\bar{\nabla}^{U}}_{\bar{D}}(\bar{D''})\rangle _U \\ 
		=&~ \langle \pi_U(\nabla^U_D D'), \pi_U(D'')\rangle_U + \langle \pi_U(D'), \pi_U(\nabla^U_D D'')\rangle_U  \\
		=&~ \rho'|_U(\langle \nabla^U_D D', D''\rangle_U + \langle D', \nabla^U_D D''\rangle_U)\\
		=&~\rho'|_U(D (\langle D', D'' \rangle_U)) \\
		=&~ \bar{D}(\langle \bar{D'}, \bar{D''}\rangle _{Y \cap U}).
	\end{align*}
\end{proof} 
Note that, the above theorem is a global counterpart of a result that appeared in \cite{VP,VP3}.
\begin{Cor}
	If an analytic subspace $Y$ of a Hermitian manifold $X$ satisfies condition (\ref{cond 1}), i.e. Remark \ref{non degenerate property} holds, then there exist an induced inner product and a Levi-civita connection on $\mathcal{T}_Y$.
\end{Cor}
 One can view these notions as generalization of the notion of Riemannian submanifolds.
 
\subsection{Relationship between Atiyah algebroid and logarithmic derivations} \label{Rel Atiyah and Logarithmic} Let $Y$ be a nonsingular divisor (submanifold) of a complex manifold $X$ of codimension $1$. Thus, there exists an open cover $\{U_{\alpha}\}_{\alpha \in \mathbb{N}}$ such that
the  divisor $Y= V(\mathcal{I})$ (vanishing set of an principal ideal sheaf $\mathcal{I}\subset \mathcal{O}_X$)  where $\mathcal{I}(U_{\alpha})=\langle f_{\alpha} \rangle$ and $X= N_Y$ (the total space of the line bundle, namely the normal bundle of $Y$ in $X$). Let $\mathcal{N}_Y$ be the associated sheaf of sections of the line bundle $N_Y \rightarrow Y$ (see \cite{AA}).

For  a (coherent) $\mathcal{O}_X$-module $\mathcal{E}$, we  forms a Lie algebroid over $(X, \mathcal{O}_X)$ consisting of the sheaf of differential operators on $\mathcal{E}$ of order $\leq 1$ with scaler symbols  \cite{BP,BRT,Ab}, i.e.
$$\mathcal{A}t(\mathcal{E}):= \{D\in \mathscr{E}nd_{\mathbb{K}_X}(\mathcal{E})~|~D(fs)=fD(s)+\sigma_D(f)s~~\mbox{ for a unique}~~  \sigma_D\in \mathcal{T}_X,~~\mbox{where}~~ f \in \mathcal{O}_X,~s\in \mathcal{E}\},$$
(thus, here $\sigma_D(f)=[D, f]_c \in \mathcal{O}_X$ for $D \in \mathcal{A}t(\mathcal{E})$ and $f \in \mathcal{O}_X$ holds) with the anchor map defined by 
\begin{center}
	$\sigma: \mathcal{A}t(\mathcal{E}) \rightarrow \mathcal{T}_X$ where $D \mapsto \sigma_D$
\end{center}
and the Lie bracket is commutator bracket.  This Lie algebroid structure is so-called Atiyah algebroid of the $\mathcal{O}_X$-module $\mathcal{E}$. 
Hence, we get the short exact sequence (s.e.s.) of Lie algebroids over $(Y, \mathcal{O}_Y)$ as follows
\vspace{-0.15cm}
$$0 \rightarrow \mathscr{E}nd_{\mathcal{O}_Y}(\mathcal{N}_Y) \hookrightarrow \mathcal{A}t(\mathcal{N}_Y)  \overset{\sigma}{\rightarrow} \mathcal{T}_Y \rightarrow 0.$$ 
Here $\mathcal{I}\subset \mathcal{O}_X$ is a locally free $\mathcal{O}_X$-module of rank $1$, we have an isomorphism of $\mathcal{O}_Y$-modules (see \cite{PT}) 
\begin{align}\label{normal Tortella}
\mathcal{A}t(\mathcal{N}_Y) \cong \mathcal{E}nd_{\mathcal{O}_Y}(\mathcal{N}_Y) \oplus \mathcal{T}_Y. 
\end{align}
More generally, if $\mathcal{I}$ is a principal ideal sheaf but not necessarily locally free then $Y$ turns out to be a hypersurface with singularities (also it is known as a principal divisor). By considering a  Hermitian inner product $\langle \cdot, \cdot \rangle : \mathcal{T}_X \times {\mathcal{T}_X} \rightarrow C^\infty_X$ on $\mathcal{T}_X$, we get the $\mathcal{O}_X$-module 
\begin{align}\label{normal sheaf}
	\mathcal{N}_{Y/X}:= \{D\in \mathcal{T}_X \mid \langle \mathcal{T}_X(- log~ Y), D \rangle \in \mathcal{I}\}.
\end{align}
It is called the sheaf of vector fields on $X$ (or derivations of $\mathcal{O}_X$) that are normal to $Y$ and the restriction of $\mathcal{N}_{Y/X}$ on $Y$ is the sheaf $\mathcal{N}_Y$.

Here, we establish an analogue of the isomorphism (\ref{normal Tortella}) in the setting of analytic spaces.
\begin{Thm} \label{log tangent-Atiyah sheaf} For a (principal) free divisor $Y$ of a Hermitian manifold $X$ (i.e. $\mathcal{T}_X(- log~ Y)$ is a locally free $\mathcal{O}_X$-module) we get the canonical isomorphism
	$$\mathcal{A}t(\mathcal{N}_{Y/X}) \cong \mathcal{T}_X(- log~ Y) \oplus \mathcal{I}~.$$
\end{Thm}
\begin{proof}
	Since Y is a principal divisor, associated to the open cover $\{U_{\alpha}\}_{\alpha \in \mathbb{N}}$ of $X$ (as mentioned in the above case) we obtain $$\mathcal{N}_{Y/X}|_{U_{\alpha}} = \langle \nabla f_{\alpha} \rangle~,$$ where $\nabla f_{\alpha} = \sum_{i=1}^{n} \partial_{x_i}(f_{\alpha})$ $\partial_{x_i}$, for all $\alpha \in \mathbb{N}$. Thus, $\nabla f_{\alpha}$ is determined by $ \{\partial_{x_i}(f_{\alpha})\}^n_{i=1}$, fully depends on $f_{\alpha}$. An element of  $\mathcal{N}_{Y/X}({U_{\alpha}})$ depends on an element of $\mathcal{I}(U_{\alpha})$ and we get an isomorphism between $\mathcal{I}$ and $\mathcal{N}_{Y/X}$ (locally defined as $gf_{\alpha} \mapsto g\nabla f_{\alpha}$ for any coefficient $g\in \mathcal{O}_X(U_{\alpha})$) of $\mathcal{O}_X$-modules. Thus, $$\mathcal{A}t(\mathcal{N}_{Y/X}) \cong \mathcal{A}t(\mathcal{I})~.$$
	
	Thus, we find local sections $D\in \mathcal{E}nd_{\mathbb{C}_X}(\mathcal{I}(U_{\alpha}))$ satisfying the property
	$$ D(g~f_{\alpha}) = g~ D(f_{\alpha}) + \sigma_{D}(g)~ f_{\alpha},$$
	for $g\in \mathcal{O}_X (U_{\alpha})$.
	 Our claim is $D \in \mathcal{A}t(\mathcal{I})$ if and only if $D \in \mathcal{T}_X(- log~ Y)(U_{\alpha})$ or $D\in \mathcal{I}(U_{\alpha})$. 
	 For any $D \in \mathcal{T}_X(- log~ Y)(U_{\alpha})$ we get $$ D(g~f_{\alpha}) = g~ D(f_{\alpha}) + D(g)~f_{\alpha} \in \mathcal{I}(U_{\alpha})~,$$ and for any $h\in\mathcal{I}(U_{\alpha})$ (note that $\mathcal{E}nd_{\mathcal{O}_X(U_{\alpha})}(\mathcal{I}(U_{\alpha})) \cong \mathcal{I}(U_{\alpha})$ canonically) we get $$h(g~f_{\alpha}) = g~hf_{\alpha} + 0~.$$
	
	 Therefore, in the global case, we get the above required isomorphism.
	\end{proof}
	\begin{Cor}
	In general, we have the s.e.s. of $\mathcal{O}_X$-modules 
		$$0\rightarrow \mathcal{I} \hookrightarrow \mathcal{A}t(\mathcal{N}_{Y/X}) \overset{\sigma}{\rightarrow} \mathcal{T}_X(- log~ Y)\rightarrow 0,$$ for any principal divisor $Y:=V(\mathcal{I})$ of a Hermitian manifold $X$.
	\end{Cor}
	\begin{Rem}
		 We have two interrelated short exact sequences
		\begin{enumerate}
			\item $0\rightarrow \mathcal{T}_X(- log~ Y) \hookrightarrow \mathcal{T}_X \overset{(\cdot)^\perp}{\rightarrow} \mathcal{N}_{Y/X} \rightarrow 0$ (of coherent $\mathcal{O}_X$-modules), and
			\item $0\rightarrow \mathcal{T}_X^0 \hookrightarrow \mathcal{T}_X(- log~ Y) \overset{\tilde{\rho}}{\rightarrow} \mathcal{T}_Y \rightarrow 0$ (of Lie algebroids over analytic spaces).
		\end{enumerate}
		Here $D^\perp := \langle D,\nabla f_{\alpha} \rangle \nabla f_{\alpha}$ for any $D\in \mathcal{T}_X(U_{\alpha})$, $[D']:=\tilde{\rho}(D')$ is the equivalence class given by the canonical map associated with  $D'\in \mathcal{T}_X(-log~ Y)(U_{\alpha})$ and $\mathcal{T}_X^0=\mathcal{I}\mathcal{T}_X$ is the space of vector fields that vanishes on $Y$ (or derivations that sends $\mathcal{O}_X$ to $\mathcal{I}$). The underlying algebraic concepts described in Section \ref{Geometry of affine variety}.
	\end{Rem}

 As in the case of complex Riemannian geometry (discussed in Section \ref{Geometry of affine variety}), consider the normal crossing divisor $Y= V(\langle xy \rangle)$  in the complex manifold $X=\mathbb{C}^2$ with standard holomorphic Riemannian metric. Then, we have $\mathcal{T}_X( - log~ Y )= \langle \{x\partial_x, y\partial_y\} \rangle$ and $\mathcal{N}_{Y/X}=\langle \{y \partial_x + x \partial_y\} \rangle$.  This is different from the usual Hermitian metric.
If we take the restrictions of these sheaves on the divisor $Y$ (that is  $\mathcal{T}_Y$ and $\mathcal{N}_Y$ respectively) and consider their stalks then both are of rank $1$ at  $Y\setminus \{(0,0)\}$. The ranks changes at the singular point $\bar{0}:=(0,0)$ where $\mathcal{T}_{Y,\bar{0}}$ is an $\mathcal{O}_{Y,\bar{0}}$-module of rank $2$ and $\mathcal{N}_{Y,\bar{0}}$ is an $\mathcal{O}_{Y,\bar{0}}$-module of rank $0$.
\section{Hermitian metrics on Lie algebroids} \label{Hermitian Lie algebroids} For a Riemannian manifold $X$, we have a symmetric positive definite bilinear form on the tangent bundle $TX$ of $X$ $($or an inner product on the $C^\infty(X)$-module of global vector fields $\mathfrak{X}(X):=\Gamma(TX))$. In the paper \cite{MB}, this notion is generalized from tangent bundle to Lie algebroids over smooth manifolds, named as Riemannian metrices on Lie algebroids. Here we consider its counter-part in complex geometry context.

\subsubsection{Hermitian Lie algebroids} Analogous to Section \ref{Hermitian metric on analytic spaces }, we define Hermitian metric on a holomorphic vector bundle and in particular on a Lie algebroid in the algebro-geometric approach, by considering positive-definite Hermitian forms on the associated sheaf of sections.

\begin{Def} \label{definition of hermitian Lie algebroids}
Let $(L, [\cdot,\cdot], \mathfrak{a})$ be a holomorphic Lie algebroid over a complex manifold $X$. It is said to be a Hermitian Lie algebroid if there exists a Hermitian metric on the underlying holomorphic vector bundle $L$.

	More generally, for a Lie algebroid $(\mathcal{L}, [\cdot,\cdot], \mathfrak{a})$  over a complex manifold $(X, \mathcal{O}_X)$, a Hermitian metric $\langle \cdot, \cdot \rangle$ on $\mathcal{L}$ is a conjugate-symmetric, positive definite, $\mathcal{O}_X$-sesquilinear map  $($or a Hermitian form$)$ given by  $($analogous notions are defined for $\mathcal{L}=\mathcal{T}_X$ in Section \ref{Hermitian metric on analytic spaces }$)$
	$$\langle \cdot, \cdot \rangle : \mathcal{L} \times {\mathcal{L}} \rightarrow C_X^\infty.$$
	
		Then $(\mathcal{L}, [\cdot,\cdot], \mathfrak{a}, \langle \cdot, \cdot \rangle)$ is said to be a Hermitian Lie algebroid over $(X, \mathcal{O}_X)$. 
\end{Def}
\begin{Rem}
    Later, we see certain compatibility conditions between $[\cdot,\cdot],~ \mathfrak{a}$ and $\langle \cdot, \cdot \rangle$ in the identity (\ref{Levi-civita-Koszul}).
\end{Rem}
\begin{Rem}
	Holomorphic Lie algebroids over a complex manifold $X$ are viewed as Lie algebroids which are locally free $\mathcal{O}_X$-modules of finite rank in the generalized context \cite{BP, BRT, AA}.
	
	A Hermitian metric $\langle \cdot, \cdot \rangle$ on a holomorphic Lie algebroid $\mathfrak{a}:L \rightarrow TX$ can be viewed as pointwise $\mathbb{C}$-sesquilinear conjugate-symmetric, positive definite maps $\langle \cdot, \cdot \rangle _x$ on the fibre $L_x$ varies smoothly on $X$ and the metric $\langle \cdot, \cdot \rangle: X \rightarrow (L \otimes \bar{L})^*$ is given by the smooth map $x \mapsto \langle \cdot, \cdot \rangle_x$. This is an analogue of the definition given in \cite{MB} for Riemannian metrics on Lie algebroids in the  complex geometry context (see Section \ref{Hermitian metric on analytic spaces }).
\end{Rem}
\begin{Exm}
	$(1)$ Let $X$ be a complex manifold. Then the tangent sheaf $\mathcal{T}_X$ of $(X, \mathcal{O}_X)$ induces a  Hermitian Lie algebroid structure  (see Section \ref{Hermitian metric on analytic spaces }). 
	\\
	$(2)$ Also, we can consider Hermitian metric on the tangent sheaf $\mathcal{L}=\mathcal{T}_Y$ of an analytic subspace $(Y, \mathcal{O}_Y)$ if it satisfies condition (\ref{cond 1}), i.e. if Remark \ref{cond 2} holds.
\end{Exm}
\subsubsection{Singular analytic foliations, leaves and invariant subspaces of a holomorphic Lie algebroid}
\begin{Def} [Singular foliations] \label{Foliation def} 
	A singular $($or regular$)$ foliation $\mathcal{F}$ on a complex manifold $X$ is a coherent $($or locally free of finite rank$)$ $(\mathbb{C}_X, \mathcal{O}_X)$-Lie-Rinehart subalgebra of $\mathfrak{X}_X$ or $\mathcal{T}_X$. It provides a generalized involutive analytic distribution on $X$ and vice-versa. $($See \cite{JV, BP, PT2}.$)$
	
	In otherwords, a singular $($or regular$)$ foliation $\mathcal{F}$ on a complex manifold $X$ is a subsheaf of $\mathcal{O}_X$-modules consists with holomorphic vector fields, which satiesfies the following properties \\
	(a) closed under the $\mathbb{C}_X$-Lie bracket $($i.e., an involutive subsheaf of $\mathfrak{X}_X$ or $\mathcal{T}_X)$,\\
	(b) locally finitely generated (free) over $\mathcal{O}_X$ $($i.e., a coherent $\mathcal{O}_X$-module$)$.
\end{Def}

	In classical differential geometric setup, a (regular) foliation over a smooth manifold $X$ represents an involutive subbundle $F \rightarrow X$ of the tangent bundle $TX \rightarrow X$, induces a Frobenious distribution. Its sheaf of sections $\Gamma_X(F)$ provides a regular foliation (as its flow) in the smooth context, as defined above. 

    \begin{Rem}A regular foliation in a complex manifold $X$ is equivalent to the sheaf of sections of a holomorphic Lie subalgebroid of the Lie algebroid $TX$ of $X$.
   \end{Rem}
\begin{Exm}\label{Foliation exm}
	For a Lie algebroid $(\mathcal{L}, \mathfrak{a},[ \cdot, \cdot])$ over a complex manifold $(X, \mathcal{O}_X)$, the image $\mathcal{F}:= \mathfrak{a}(\mathcal{L}) \subset \mathcal{T}_X$ of the anchor map forms a singular foliation.
\end{Exm} 
\begin{Def} [Leaves or Orbits] \label{Leaves}
	An integral submanifold of a foliation $\mathcal{F}$ $($in both real smooth or complex analytic cases$)$ is an immersed submanifold $Z \subset X$ $($where $X$ is smooth or complex manifold respectively$)$ with the property that
	for all $p \in Z$, $$T_p Z= Span\{D(p) \in T_p X \mid  D \in \mathcal{F}\}.$$ A leaf (or  orbit) is a maximal connected integral submanifold.
\end{Def}

\begin{Thm}\label{Stefan-Sussman-Nagano}
	$(1)$ $($Stefan-Sussman's theorem$)$ If $\mathcal{F}$ is a singular foliation on a real smooth manifold $X$ then there is a leaf through every point in $X$ (see \cite{RF}).\\
	$(2)$ $($Nagano's theorem$)$ If $\mathcal{F}$ is a singular analytic distribution (foliation) on a complex manifold $X$ then $X$ is partitioned into leaves (see \cite{BP}).
\end{Thm}
\begin{Rem}\label{L-invariant}
	In general, it is useful to consider analytic subspaces $Y \subset X$ that are union of leaves of the foliation $\mathfrak{a}(\mathcal{L})$ for a Lie algebroid $\mathcal{L}$. These are $\mathcal{L}$-invariant subspaces,  appear to study all of the vector fields coming from  $\mathcal{L}$ that are tangent to $Y$. More explicitly, if a closed complex analytic subspace $Y:=V(\mathcal{I})$ of some ideal sheaf $\mathcal{I}$ satisfies $\mathfrak{a}(\mathcal{L})(\mathcal{I}) \subset \mathcal{I}$, then $Y$ is called an \emph{$\mathcal{L}$-invariant subspace} (see \cite{BP}).
\end{Rem}
	Frobenius distributions and Stefan-Sussmann distributions are represents integrable regular and singular foliations respectively over a smooth manifold $(X, C_X^\infty)$, where $C_X^\infty$ is the sheaf of real valued $C^{\infty}$-functions on $X$. These corresponds with Lie algebroids over $(X, C_X^\infty)$, namely generalized involutive subsheaf of $C_X^\infty$-modules of the tangent sheaf $\mathcal{T}_X$. 
    In complex geometry its analogue provides a particular kind of holomorphic foliations (or analytic distributions), which forms Lie algebroids over a complex manifold $(X,\mathcal{O}_X)$. 

\subsection{Characteristic foliation with induced inner product} \label{Characteristic foliation} Let $(L, [\cdot,\cdot], \mathfrak{a})$ be a holomorphic Lie algebroid over a complex manifold $X$ \cite{BP}. Thus its image of the anchor map $\mathfrak{a}: L \rightarrow TX$ defines a singular foliation of $X$. Consider the sheaf of holomorphic sections $\Gamma_L$ of the map $p:L \rightarrow X$, denote it by $\mathcal{L}$. 
 Since the sheaf homomorphism $\mathfrak{a}$ is an $\mathcal{O}_X$-linear and $\mathbb{C}_X$-Lie algebra homomorphism, for any two section $\tilde{D}_1, \tilde{D}_2 \in \mathfrak{a}(\mathcal{L})$ and a section $f \in \mathcal{O}_X$ we obtain the sections of $\mathfrak{a}(\mathcal{L})$
 $$f\tilde{D}_1 = \mathfrak{a}(f D_1)~\text{and}~[\tilde{D}_1, \tilde{D}_2] =\mathfrak{a}([D_1,D_2]),$$
 where $D_1, D_2$ are some preimage of $\tilde{D}_1, \tilde{D}_2$ respectively. 
Thus the image $\mathscr{I}m(\mathfrak{a})$ or $\mathfrak{a}(\mathcal{L})$ of the anchor map forms a coherent involutive subsheaf of the tangent sheaf $\mathcal{T}_X$ defined by the sheafification of the presheaf $U \mapsto Im(\mathfrak{a}(U))$, called characteristic foliation \cite{MB} or orbit foliation \cite{RF} of $L$. In the complex geometry context, we use Nagano's theorem on the integrability of singular analytic distributions (Theorem \ref{Stefan-Sussman-Nagano}). 

We denote by $L_Z$ the restriction of the holomorphic Lie algebroid $L$ to a leaf $Z$.  One can deduce easily that the bracket $[\cdot,\cdot]$ induces a bracket on the space of holomorphic sections of $p_Z : L_Z \rightarrow Z$ where $p : L \rightarrow X$ is the projection and $p|_Z:=p_Z$, and hence a transitive Lie algebroid structure.

Let $\langle \cdot, \cdot\rangle$ be a Hermitian metric on the Lie algebroid $\mathfrak{a} :L \rightarrow TX$. Then for any leaf $Z$ of the characteristic foliation and for any $x \in Z$, the $\mathbb{C}$-inner product space  $$L_x = \mathcal{G}^0_x \oplus \mathcal{G}^T_x,$$ where $L_x=p^{-1}(x)$ and $\mathcal{G}^T_x$ is the orthogonal complement to $\mathcal{G}^0_x = ker(\mathfrak{a}_x)$ with respect to $\langle \cdot, \cdot \rangle_x$, where $\mathfrak{a}_x$ is the fiberwise map at $x$. Now the restriction of the map $\mathfrak{a}_x$ to $\mathcal{G}^T_x$ is an isomorphism into $T_xZ$ and hence induces a scalar product on $T_xZ$, is given by $\langle v, v' \rangle_{T_xZ}:= \langle b, b'\rangle_x,$ where $b,b' \in \mathcal{G}^T_x$ with $\mathfrak{a}_x(b) = v$ and $\mathfrak{a}_x(b')= v'$, which varies smoothly on $Z$. Thus $\langle \cdot, \cdot \rangle$ induces a Hermitian metric $\langle \cdot, \cdot \rangle_Z$ on $Z$.
Fix a leaf $Z$ and consider $p_Z: L_Z \rightarrow Z$. We have $$L_Z = \mathcal{G}^0_Z \oplus \mathcal{G}^T_Z$$ and we call the elements of $\Gamma(\mathcal{G}^0_Z)$ vertical sections and the elements of $\Gamma(\mathcal{G}^T_Z)$ horizontal sections. Thus, we get a short exact sequence of Hermitian (holomorphic) Lie algebroids
\begin{center}
	$0 \rightarrow \mathcal{G}^0_Z \rightarrow L_Z \rightarrow TZ \rightarrow 0$
\end{center}
is formally identical to a Hermitian submersion.

\subsubsection{Induced metric on invariant subspaces}
The characteristic foliation of a holomorphic Lie algebroid $\mathfrak{a}: L \rightarrow TX$ is associated with a system of first order partial differential equations, provides a family of leaves as integral solutions and forms $L$-invariant subspaces (see Remark \ref{L-invariant}). For a $\Gamma_X(L)=:\mathcal{L}$-invariant subspace $Y$ in $X$, the foliation $\mathscr{I}m(\tilde{\mathfrak{a}})$ is a Lie subalgebroid of $ \mathcal{T}_X(-log~ Y)$, where $\tilde{\mathfrak{a}}=\Gamma_X(\mathfrak{a})$. Our aim is to show that there is a Hermitian metric on $Y$ (or Hermitian inner product on $\mathcal{T}_Y)$ induces from the Hermitian metric on the Lie algebroid.
First notice that as $\mathcal{O}_X$-modules we obtain $$\mathcal{L}= \mathscr{K}er(\mathfrak{a}) \oplus (\mathscr{K}er(\mathfrak{a}))^{\perp},$$  where the $\mathcal{O}_X$-module $(\mathscr{K}er(\mathfrak{a}))^{\perp}$ is the sheafification of the presheaf
$$U \mapsto \{D\in \mathcal{L}(U)~|~ \langle Ker(\mathfrak{a}(U)), D \rangle=0 \}.$$
Thus, for any section $D$ of $\mathcal{L}$ (i.e. $D \in \mathcal{L}$), we have $D_1 \in \mathscr{K}er(\mathfrak{a})$ and $D_2 \in (\mathscr{K}er(\mathfrak{a}))^{\perp}$ with $D=D_1+D_2$ or $D=(D_1, D_2)$. In this situation, we get $\mathfrak{a}(D)= \mathfrak{a}(D_2)$, implies $\mathfrak{a}(\mathcal{L})= \mathfrak{a} ((\mathscr{K}er(\mathfrak{a}))^{\perp})$.
\begin{Thm} \label{metric on invariant subspaces}
	An $\mathcal{L}$-invariant subspace $Y \subset X$ of a Hermitian $($holomorphic$)$ Lie algebroid $\mathfrak{a} :\mathcal{L} \rightarrow \mathcal{T}_X$ canonically induces a Hermitian metric on $Y$ if it satisfies $\mathfrak{a}(\mathcal{L}) =  \mathcal{T}_X(-log~ Y)$.
\end{Thm}
\begin{proof}
	Let $D, D' \in \mathfrak{a}(\mathcal{L}) =  \mathcal{T}_X(-log~ Y)$, then we define an inner product on $\mathfrak{a}(\mathcal{L})$ as
	$$\langle D, D' \rangle _{\mathfrak{a}(\mathcal{L})} := \langle \tilde{D}, \tilde{D'} \rangle$$ where $\tilde{D}, \tilde{D'} \in (\mathscr{K}er(\mathfrak{a}))^{\perp} \subset \mathcal{L}$ such that $ \mathfrak{a}(\tilde{D})= D,~ \mathfrak{a}(\tilde{D'})= D'$ (without loss of generality). To show that the map is well defined, consider $\bar{D}, \tilde{D} \in \mathcal{L}$ such that $ \mathfrak{a}(\tilde{D})= D= \mathfrak{a}(\bar{D})$ and $\bar{D}' \in (\mathscr{K}er(\mathfrak{a}))^{\perp}$ with $\mathfrak{a}(\bar{D}')=D'$. Thus, $\mathfrak{a}(\bar{D}-\tilde{D})=0$, implies $\bar{D}-\tilde{D} \in \mathscr{K}er(\mathfrak{a})$, and hence $\langle \bar{D}-\tilde{D}, \bar{D}'\rangle=0$, implies $\langle \bar{D}, \bar{D}'\rangle = \langle \tilde{D}, \bar{D}'\rangle.$

	Thus, the surjective $\mathcal{O}_X$-algebra homomorphism $\rho': C_X^\infty \rightarrow C_X^\infty/{\mathcal{I}}$ (see Notation \ref{Notation}) and the surjective Lie algebroid homomorphism $\tilde{\rho}: \mathcal{T}_X(-log Y) \rightarrow \mathcal{T}_Y$ canonically induces  a Hermitian inner product on $\mathcal{T}_Y$ as 
	$$\langle D, D' \rangle _{Y} := \rho'( \langle \tilde{D}, \tilde{D'} \rangle_{\mathfrak{a}(\mathcal{L})})~,$$
	where $\tilde{D}, \tilde{D'} \in \mathcal{T}_X(-log~Y)$ such that $ \tilde{\rho}(\tilde{D})= D,~ \tilde{\rho}(\tilde{D'})= D'$ (see Theorem \ref{induced inner product} for details).
\end{proof}
 \begin{Cor} We obtain a short exact sequence of Hermitian Lie algebroids over $(X, \mathcal{O}_X)$ as
 	\begin{center}
 		$0 \rightarrow \mathscr{K}er(\mathfrak{a})\hookrightarrow \mathcal{L}  \overset{\mathfrak{a}}{\rightarrow} \mathscr{I}m(\mathfrak{a}) \rightarrow 0$~.
 	\end{center}
 \end{Cor}

\begin{Rem}
Let $Y=V(\mathcal{I})$ be a hypersurface with isolated singularities of $X$. Consider the stalks $\mathcal{T}_{X,x}(-log~ Y)$ of the Lie algebroid $\mathcal{T}_X(-log~ Y)$ of sheaf of holomorphic vector fields on $X$ that are tangent to $Y$ (or derivations of $\mathcal{O}_X$ that preserves the ideal sheaf $\mathcal{I}$). This Lie algebroid associated with the characteristic foliation and the corresponding stalks change its rank if it is not regular. Thus leaves (or orbits) of this foliation are of different dimensions and forms $L$-invariant subspace $Y$ as a singular analytic space. 	
\end{Rem}

\subsubsection{Levi-Civita Connection on a Hermitian Lie algebroid} We extend the notion of covariant connection for a hermitian manifold
(see Section \ref{connection}), not necessarily a holomorphic connection but compatible with the metric, in the context of holomorphic Lie algebroids. Here, we define an analogue of Levi-civita connection for a Hermitian Lie algebroid as follows.

Let $\mathcal{L}$ be a holomorphic Lie algebroid over a complex manifold $(X, \mathcal{O}_X)$. Since $\mathcal{L}$ is a locally free $\mathcal{O}_X$-module of some finite rank, say $n$, there exist an open cover $\mathcal{U}:=\{U_i\}_{i\in \mathbb{N}}$ of $X$ such that
$\mathcal{L}|_{U_i}$ is a free $\mathcal{O}_X|_{U_i}$-module of rank $n$, for each $i \in \mathbb{N}$.
Consider a $C_X^\infty$-module homomorphism 
\begin{align} \label{complexified smooth connection}
	\nabla: {C^{\infty}_X} \otimes_{\mathcal{O}_X} \mathcal{L} \rightarrow {C^{\infty}_X} \otimes_{\mathcal{O}_X} \mathcal{A}t(\mathcal{L}) 
\end{align}
such that there exists holomorphic Lie algebroid connection $\nabla^i$ (see \cite{PT}), i.e. a Lie algebroid homomorphism $$\nabla^{i}: \mathcal{L}|_{U_i} \rightarrow \mathcal{A}t(\mathcal{L}|_{U_i}),$$ for each $i \in \mathbb{N}$,  determines $\nabla$ by globalizing with a partition of unity of $\mathcal{U}$. 
This is different from usual $\mathcal{L}$-connection on $\mathcal{L}$. A holomorphic Lie algebroid connection always produce such a connection.

For example, consider the $\mathcal{L}$ connection $\nabla := \sum_{i} f_i ~ \nabla^{i}$, i.e. $\nabla_D(s)= \sum_{i} f_i ~ \nabla^{i}_D(s)$, for $D, s \in \mathcal{L}(U_i)$ where $\{f_i\}_{i\in \mathbb{N}}$ is a smooth partition of unity for the open cover of $X$. 
Consider a basis $\{D^i_1, \dots, D^i_n\}$  of the $\mathcal{O}_X(U_i)$-module $\mathcal{L}(U_i)$ for each $U_i \in \mathcal{U}$, and thus the map $\nabla^{U_i}_D$ is defined by 
$$\nabla^{i}_D( \sum^n_{j=1} g_j~ D^i_j) = \sum^n_{j=1} \mathfrak{a}(D)(g_j)~ D^i_j.$$

 An analogue of the Levi–Civita connection for Hermitian metrics on Lie algebroids is an important ingredient (see \cite{MB}). If $\langle \cdot,\cdot \rangle$ is a Hermitian metric on a Lie algebroid $\mathfrak{a}: L \rightarrow TX$, then the Koszul formula 
   \begin{equation}\label{Levi-civita-Koszul}
\begin{aligned} 
    2\langle \tilde{\nabla}_D D', D'' \rangle &= \mathfrak{a}(D)(\langle D', D'' \rangle) + \mathfrak{a}(D')(\langle D, D'' \rangle) - \mathfrak{a}(D'')(\langle D, D' \rangle) \\
    &\quad - \big( \langle D,[D', D''] \rangle +  \langle D',[D, D''] \rangle - \langle D'',[D, D'] \rangle \big).
\end{aligned}
\end{equation}
holds for sections $D,D',D''$ are  in $\mathcal{L}:=\Gamma_X(L)$,
defines a linear $\mathcal{L}$-connection on $\mathcal{L}$ in the sense of $(\ref{complexified smooth connection})$, which is characterized by the following two properties  (on each space of sections): 
\begin{itemize}
    \item $\tilde{\nabla}$ is compatible with the metric: $\mathfrak{a}(D)(\langle D', D''\rangle) = \langle \tilde{\nabla}_D D', D'' \rangle + \langle D', \tilde{\nabla}_D D'' \rangle$,
    \item $\tilde{\nabla}$ is torsion free or symmetric: $\tilde{\nabla}_D D' - \tilde{\nabla}_{D'} D = [D, D']$.
\end{itemize}
We call $\tilde{\nabla}$ is a Levi-Civita $L$ -connection associated to the Hermitian metric $\langle \cdot, \cdot \rangle$.

If $\nabla$ be a $TX$-connection on $L$ then there are two obvious linear $L$-connection
\begin{center}
	$\tilde{\nabla}^0_D D' = \nabla_{\mathit{a}(D)} D'$ and  $\tilde{\nabla}^1_D D' = \nabla_{\mathit{a}(D')} D + [D, D']$,
\end{center}
which are same for the standard case of the Lie algebroid $TX$ with the Levi-Civita connection.
\begin{Prop}
	An $\mathcal{L}$-invariant subspace $Y \subset X$ of a Hermitian $($holomorphic$)$ Lie algebroid $\mathcal{L}$ over $(X, \mathcal{O}_X)$, canonically induces a $\mathcal{T}_Y$-Levi-civita connection with the induced metric on $Y$.
\end{Prop}

\begin{Exm} \label{Foliation with inner product}
	\normalfont{In the classical set up, the cotangent bundle over a Poisson manifold has a natural Lie algebroid structure \cite{KM}. The associated sheaf of sections provides a Lie algebroid in the general context. Here, we consider the case associated with an analytic space, namely the nilpotent cone (see \cite{BP,AA}).
	
	Let us consider the complex manifold $\mathbb{C}^3$ with the Poisson algebra structure on $(\mathcal{O}_{\mathbb{C}^3},\{\cdot,\cdot\})$, denoted by $X$ \cite{BP,AA}. This structure is induced by the standard Lie algebra $(\mathfrak{sl}_2(\mathbb{C}),[\cdot,\cdot]_c)$
	One may describe the nilpotent cone through a generalized involutive distribution generated by the Hamiltonian vector fields in the Poisson manifold $X$.
	The (sheaf of) Poisson algebra $(\mathcal{O}_X, \{\cdot,\cdot \})$ defined on coordinate functions $x, y, z$ (over any open subset of $\mathbb{C}^3$) as follows. $$\{x,y\} = 2y, \{x,z\} = -2z, \{y,z\} = x.$$
	Then the associated Hamiltonian vector fields over $X$  are given by 
    \begin{align} \label{Hamiltonian vector fields} 
       D_x = 2y\partial_y - 2z \partial_z ,~  D_y = -2y\partial_x + x \partial_z ,~ D_z =  2z\partial_x - x \partial_y. 
    \end{align} 
	We consider the system of first order homogeneous  partial differential equations
	$$D_x(f) = D_y(f) = D_z(f) = 0.$$ 
	Then we find solution of the above system is $f(x,y,z) = x^2 +4yz$. It describes a family of level sets $f^{-1}(c)$, is  parametrized by $c \in \mathbb{C}$. For each nonzero values of $c$, we get  a $2$-dimensional submanifold.  Now for $c=0$, $Y:= f^{-1}(0)$ is the vanishing set of the ideal sheaf $\mathcal{I}:=\langle x^2 + 4yz \rangle \subset \mathcal{O}_{X}$, is hypersurface with singularity at the origin, known as the nilpotent cone.
	
	It yields  a Lie algebroid structure on the cotangent sheaf 
    $\Omega^1_{X}$, where the anchor map 
    $$\mathfrak{a}:\Omega^1_{X}\rightarrow \mathcal{T}_{X}~\text{is locally given by}~f dg \mapsto f\{g,\cdot \},~\text{for any}~f,g \in \mathcal{O}_{X}$$ and the Lie algebra structure is locally given by (using Lie derivatives)
    $$[df,dg]=\mathfrak{L}_{\mathfrak{a}(df)}(dg)-\mathfrak{L}_{\mathfrak{a}(dg)}(df)-d(\{f,g\}),~\text{for any}~f,g \in \mathcal{O}_{X}.$$ 
	The sheaf of sections of the cotangent bundle $T^*X$ (forms a holomorphic Lie algebroid) is $\Omega^1_{X}$, provides the Lie algebroid  $\mathscr{I}m(\mathfrak{a})$ in the general set up.  It is a generalized involutive subsheaf of the tangent sheaf $\mathcal{T}_{X}$. 
	As a Lie algebroid, the $\Omega^1_{X}$-invariant subspaces \cite{BP} are the hypersurface singularities, one of them is the nilpotent cone $Y$  associated with the principal ideal sheaf $\mathcal{I}= \langle x^2+ 4yz \rangle$. 
	The associated integrable system of first order homogeneous  PDE's, consists of the Hamiltonian vector fields, is given as follows
	\begin{center}
		$\mathfrak{a}(dx)=\{x, \cdot \}=D_x,~
		\mathfrak{a}(dy)=\{y, \cdot \}=D_y,~
		\mathfrak{a}(dz)=\{z, \cdot \}=D_z,$
	\end{center}
	corresponds to the characteristic foliation $Im (\mathfrak{a}) = \mathcal{T}_{X}(-log~ Y)$. 
	
	Consider the standard Hermitian metric (see \cite{VP}) on the Lie algebroid $\Omega^1_{X}$, given as 
	$$\langle dx, dx \rangle =\langle dy, dy \rangle=\langle dz, dz \rangle=1;~ \langle dx, dy \rangle=\langle dx, dz \rangle=\langle dy, dz \rangle=0.$$
	This provides an $\mathcal{O}_{X}$-sesquilinear form on the characteristic foliation $\mathcal{T}_{X}(-log~ Y)$, given by
	$$\langle D_x, D_x \rangle =\langle D_y, D_y \rangle=\langle D_z, D_z \rangle=1;~ \langle D_x, D_y \rangle=\langle D_x, D_z \rangle=\langle D_y, D_z \rangle=0.$$
	Since the $\mathcal{O}_{X}$-module $\mathcal{T}_{X}(-log~ Y)$ is generated by the derivations $D_x,D_y,D_z$, thus any section of it is of the form $D=f D_x+ g D_y + h D_z$ for some $f,g,h \in \mathcal{O}_{\mathbb{C}^3}$. To show that the induced inner product is positive definite or non degenerate, notice that 
    $$\langle D, D_x \rangle= \langle D, D_y \rangle=\langle D, D_z \rangle=0, ~\text{implies}~f=g=h=0, ~\text{i.e.}~ D=0$$ and 
    $$\langle D, D \rangle= |f|^2+|g|^2+|h|^2 =0, ~\text{implies}~D=0.$$  Hence, $\mathcal{T}_{X}(-log~ Y)$ is a Hermitian Lie algebroid over the complex manifold $(X,\mathcal{O}_{X})$ with the induced Hermitian inner product.
	
	Now, we consider the standard inner product $\langle \cdot, \cdot \rangle_{\mathbb{C}^3}$ of $\mathcal{T}_{\mathbb{C}^3}$ (described in Remark \ref{standard inner product}) and take its restriction on  $\mathcal{T}_{\mathbb{C}^3}(-log~ Y)$, to show that it  differs from the above described inner product. 
	Thus, from the equations (\ref{Hamiltonian vector fields}) we get the followings:
	$$\langle D_x, D_x \rangle_{\mathbb{C}^3}= 4(|y|^2+  |z|^2), ~\langle D_y, D_y \rangle_{\mathbb{C}^3}= |x|^2+ 4 |y|^2, ~\langle D_z, D_z \rangle_{\mathbb{C}^3}= |x|^2+4 |z|^2,$$ 
     $$\langle D_x, D_y \rangle_{\mathbb{C}^3} =-2z \bar{x}, ~\langle D_x, D_z \rangle_{\mathbb{C}^3}=-2y \bar{x}, ~\langle D_y, D_z \rangle_{\mathbb{C}^3} = -4y \bar{z}.$$
	
	Therefore, the induced inner product on the characteristic foliation  $Im (\mathfrak{a}) = \mathcal{T}_X(-log~ Y) $ differs from the standard inner product on  $\mathcal{T}_X(-log~ Y)\subset \mathcal{T}_X$.}
\end{Exm}

\section{The analytic de Rham cohomology for a Characteristic Foliation} \label{The analytic de Rham cohomology of spaces associated with Characteristic Foliation} 
Here, we examine cohomological aspects of the characteristic foliation introduced in Section \ref{Characteristic foliation}, which is associated with a (holomorphic) Lie algebroid over a Hermitian manifold in the framework of complex geometry. We study the cohomology theory of the induced leaf space, individual leaves, and certain unions of leaves. This serves as our model for singular analytic foliations on a Hermitian manifold. In the case of a regular foliation on a Riemannian manifold, a notion of de Rham cohomology of the leaf space \cite{GS} and the cohomology of individual leaves \cite{JLH} have been explored. We extend this study to the holomorphic setting, developing analogous cohomological structures.

\subsection{Cohomology for Leaf space } Let $X$ be a Hermitian manifold and $\mathcal{F}$ a  foliation on $X$. 
Let $\Omega^k_{\mathcal{F}}$ denote the sheaf of holomorphic differential $k$-forms on $X$ $(k \in \mathbb{N} \cup \{0\})$ which are locally satisfying the following identities (using contraction map $i$, exterior differentiation $d$ and Lie derivative $\mathfrak{L}$)
\begin{center}
	$i_{D} \omega = 0 = i_{D} d{\omega}$,\\
	\hspace{-.6 cm}	or,~ ~$i_{D} \omega = 0 = \mathfrak{L}_{D} \omega$, 
\end{center}
for a section $D$ of $\mathcal{F}$ and a section $\omega$ of $\Omega^k_X$. Thus the sheaf $\Omega^k_{\mathcal{F}}$ consists of holomorphic differential $k$-forms which are invariant under the flow (induced by $\mathcal{F}$).
Now, if we consider  the characteristic foliation $\mathcal{F}=\mathfrak{a}(\mathcal{L})$ of a holomorphic Lie algebroid $\mathfrak{a}:\mathcal{L}\rightarrow \mathcal{T}_X$ then using the Nagano's integrability theorem (see Theorem \ref{Stefan-Sussman-Nagano}) we get a family of analytic spaces which are maximal $\mathcal{L}$-invariant subspaces associated with the foliation $\mathcal{F}$, denote one of them as $Y$. Then the $\mathcal{O}_X$-module $$\Omega^1_{\mathcal{F}} := (\mathcal{T}_X/{\mathcal{F}})^*=\mathscr{H}om_{\mathcal{O}_X}(\mathcal{T}_X/{\mathcal{F}},~\mathcal{O}_X)$$ is the conormal sheaf of $\mathcal{F}$. In the special case when $\mathfrak{a}(\mathcal{L}) = \mathcal{T}_X(-logY)$ (for e.g. see Example \ref{Foliation with inner product}), 
the normal sheaf $\mathcal{N}_{Y/X}$ (defined in (\ref{normal sheaf})) of $Y$ in $X$ is equivalent to $\mathcal{T}_X/\mathcal{F}$ for an analytic subspace $Y$ of $X$ as described above. The associated exterior algebra  over $\mathcal{O}_X$ forms a cochain complex of $\mathcal{O}_X$-modules $$\Omega^\bullet_{\mathcal{F}}:= (\wedge^\bullet_{\mathcal{O}_X}\Omega_{\mathcal{F}}^1,~ d_{\mathcal{F}}),$$
where $d_{\mathcal{F}}$ is the canonical differential on $\wedge^\bullet_{\mathcal{O}_X}\Omega_{\mathcal{F}}^1$ induced from the differential $d_X$ of the analytic de Rham complex $\Omega^\bullet_X:=(\wedge^\bullet_{\mathcal{O}_X}\Omega_{X}^1,~ d_X)$ of $X$ (see \cite{BRT, Ab}). The exterior differentiation $d_{\mathcal{F}}$ leaves $\Omega^\bullet_{\mathcal{F}}$ invariant, and we denote  the hypercohomology ring of the cochain complex $\Omega^\bullet_{\mathcal{F}}$ by $\mathbb{H}^\bullet_{\mathcal{F}}(X)$.

\begin{Rem}The leaf space \( X/{\mathcal{F}} \) is a topological space obtained by collapsing each leaf of \( \mathcal{F} \) to a single point. The cochain complex of \( \mathcal{O}_X \)-modules \( \Omega^\bullet_{\mathcal{F}} \) can be viewed, in a certain sense, as the sheaf of algebras of holomorphic differential forms along the leaves of \( \mathcal{F} \). The associated cohomology ring \( \mathbb{H}^\bullet_{\mathcal{F}}(X) \) can be interpreted as an analytic counterpart to the foliated de Rham cohomology of \( X \).  
\end{Rem}
\subsection{Cohomology of leaves (or orbits) } \vspace{-.1 cm}Consider the normal sheaf to a foliation $\mathcal{F}$ in the Hermitian Lie algebroid $\mathcal{T}_X$, denote it by $\nu$ and thus $\nu = \mathcal{T}_X/{\mathcal{F}}$. If $\mathcal{F}$ is regular then there is a canonical connection $\nabla$ on $\nu$ which is flat along $\mathcal{F}$ (sheafifiying the smooth Bott connection \cite{JLH} and considering its analogue in holomorphic context), i.e. $\nu$ is a $\mathcal{T}_X$-module and the associated cotangent sheaf to the foliation forms a  cochain complex of sheaves $$\Omega^\bullet _{\mathcal{F}}(\nu)=(\wedge^\bullet_{\mathcal{O}_X}(\mathcal{F}^*)\otimes_{\mathcal{O}_X} \nu, ~\tilde{d})$$ where $\wedge^\bullet_{\mathcal{O}_X}(\mathcal{F}^*)\otimes_{\mathcal{O}_X} \nu$
is the sheafification of the presheaf $$U\mapsto Hom_{\mathcal{O}_X(U)}(\wedge^{\bullet}_{\mathcal{O}_X(U)}\mathcal{F}(U),~ \mathcal{\nu}(U))= \wedge^\bullet_{\mathcal{O}_X(U)}((\mathcal{F}(U))^* \otimes_{\mathcal{O}_X(U)} \nu(U)$$
(since here $\mathcal{F}$ and $\nu$ both are locally free $\mathcal{O}_X$-module of finite rank) and the differential $\tilde{d}$ is the Chevalley-Eilenberg-de Rham differential \cite{Ab} for the Lie algebroid $\mathcal{F} \subset \mathcal{T}_X$ with coefficient in the $\mathcal{F}$-module $(\nu, \nabla|_{\mathcal{F}})$.

If  the characteristic foliation of a  holomorphic Lie algebroid $\mathcal{L}$ is of the form  $\mathcal{T}_X(-log~Y)$ for some $\mathcal{L}$-invariant subspace $Y \subset X$, then the cochain complex is the Lie algebroid complex of $\mathcal{T}_X(-log~Y)$ with coefficients in the $\mathcal{O}_X$-module $\nu = \mathcal{N}_{Y/X}$. In particular, for a free  divisor $Y$ (i.e. $\mathcal{T}_X(-log~Y)$ is a locally free $\mathcal{O}_X$-module), this complex is consists of sheaf of logarithmic differential forms with coefficients in $\mathcal{N}_{Y/X}$
\begin{align} \label{log differntial forms}
  \Omega^\bullet_X(log ~Y) \otimes_{\mathcal{O}_X} \mathcal{N}_{Y/X},  
\end{align}
its (global) sections are meromorphic forms on $X$ with simple poles along the divisor $Y$  (see \cite{PT2, Ab}).

When $\mathcal{F}$ is a regular foliation, the restriction $\mathcal{F}|_Z$ and $\nu|_Z$ on a leaf $Z$ are sheaf of sections of the tangent bundle and the normal bundle of $Z$ respectively. By extending the analytic de Rham theorem,  locally using complex geometric analogue of the case as described in \cite{JLH}, we get the following result. 
\begin{Prop} \label{cohomology of a leaf}
	The hypercohomology of  the cochain complex of $\mathcal{O}_X$-modules $(\wedge^\bullet_{\mathcal{O}_Z}(\mathcal{T}^*_Z)\otimes_{\mathcal{O}_Z} \nu|_Z, \tilde{d})$ is isomorphic $($as graded $\mathbb{C}$-vector spaces$)$ to the singular cohomology $H^\bullet (Z, \mathbb{C}^q)$ where $q = dim$ $\nu$.
\end{Prop}
More generally, for a Lie algebroid $\mathcal{L}$ with coefficient in $(\mathcal{E}, \nabla)$, we have the \emph{Chevalley-Eilenberg-de Rham differential} $d_{\mathcal{L}}: \wedge^\bullet_{\mathcal{O}_X}\mathcal{L}^* \otimes_{\mathcal{O}_X} \mathcal{E} \rightarrow \wedge^{\bullet +1}_{\mathcal{O}_X}\mathcal{L}^* \otimes_{\mathcal{O}_X} \mathcal{E}$ (see \cite{BRT, PT2}), recovers the earlier cases,  is given by
\begin{align*}
	\begin{split}
		d_{\mathcal{L}}(\omega)(D_1\wedge \cdots \wedge D_{k+1}) & = \sum^{k+1}_{i=1}(-1)^{i+1}~\nabla_{D_i}(\omega(D_1\wedge \cdots \wedge \hat{D_i} \wedge \cdots \wedge D_{k+1}))\\
		&+ \sum_{i< j}(-1)^{i+j}~\omega([D_i, D_j]\wedge D_1\wedge \cdots \wedge \hat{D_i}\wedge \cdots \wedge \hat{D_j}\wedge \cdots \wedge D_{k+1}),
	\end{split}
\end{align*}
where $D_1, \dots, D_{k+1} \in \mathcal{L}$ and
$\omega \in \wedge^{k}_{\mathcal{O}_X} \mathcal{L}^*$, 
and $\nabla: \mathcal{L} \rightarrow \mathcal{A}t(\mathcal{E})$ is the flat $\mathcal{L}$-connection on $\mathcal{E}$ (see Section \ref{Rel Atiyah and Logarithmic}).

\section{Equivariant $($holomorphic$)$ Lie algebroid cohomology} \label{Equivarient Lie algebroid cohomology}
In \cite{EE}, the Lie-Rinehart cohomology of quotient singularities by finite groups was studied in the context of affine algebraic geometry. Moreover, in \cite{BR}, the equivariant cohomology of smooth Lie algebroids was examined. In this section, we develop an analogous framework in both the smooth and analytic settings.

\subsection{Algebraic structures induced by a Lie group action} Let $X$ be a smooth manifold carrying an action of a Lie group $G$. Then the Lie group action can be expressed as a group homomorphism
\begin{center}
	$\phi: G \rightarrow Aut_{\mathbb{R}}(C^\infty (X))$ defined by $g\mapsto \phi_g$, 
\end{center}
where for all $f \in C^\infty (X)$, $\phi_g(f) = g\cdot f$  and $(g\cdot f)(x)= f(g^{-1}\cdot x)$ for all $x \in X$. Then $\phi$ naturally induces a group action (known as the conjugate action) of $G$ on $Der_{\mathbb{R}}(C^\infty(X))$, given by the action
\begin{center}
	$g\cdot D = gDg^{-1}$ $~$ for all $g \in G$ and $D \in Der_{\mathbb{R}}(C^\infty(X))$
\end{center}
where $gDg^{-1}(f)=  g\cdot(D(g^{-1}\cdot f))$, for all $f \in C^\infty (X)$.
Therefore,  we obtain the identity
\begin{center}
	$g\cdot(fD) = (g\cdot f)(g\cdot D)$ for all $g \in G$, $f \in C^\infty (X)$ and $ D \in Der_{\mathbb{R}}(C^\infty (X))$.
\end{center}

\subsubsection{$C^\infty (X)-G$ module structures}
We call a $C^\infty (X)$-module $M$ has a $C^\infty (X)-G$ module structure if there is a group action $G\mapsto Aut_{\mathbb{R}}(M)$ which is compatible with the $C^\infty(X)$-module structure, i.e. a group action such that $g\cdot (fm) = (g\cdot f)(g\cdot m)$ for all $g\in G, f\in C^\infty(X), m\in M$. Hence, $Der_{\mathbb{R}}(C^\infty (X))\cong \mathfrak{X}(X)$ has a natural $C^\infty (X)-G$ module structure induced by the map $\phi$, where $\mathfrak{X}(X)$ is the space of (smooth) vector fields on $X$.
For $C^\infty (X)-G$ modules $M$ and $N$, consider  the $C^\infty(X)$-modules $Hom_{C^\infty (X)}(M,N)$ and on $M\wedge_{C^\infty(X)} N$ with the (conjugate) $G$-actions
\begin{center}
	$g\cdot \omega = g\omega g^{-1}$ where ~ $ g\omega g^{-1}(m)= g\cdot(\omega(g^{-1}\cdot m))$, and $g\cdot (m\wedge n) = (g\cdot m)\wedge (g\cdot n)$,
\end{center}
for all $g\in G, \omega \in Hom_{C^\infty (X)}(M,N)$, $m\in M$ and $n\in N$. This provides a canonical $C^\infty (X)-G$ module structures on the $C^\infty(X)$-modules. Thus, the space of differential $k$-forms is defined by 
$$\Omega^k(X):= Hom_{C^\infty (X)}(\wedge^k_{C^\infty(X)}\mathfrak{X}(X), C^\infty(X))$$
has a  natural $C^\infty (X)-G$ module structure induced by the map $\phi$. 
\subsubsection{Group of diffeomorphisms}
In differential geometry, the notion of a smooth manifold $X$ carries an  action of a Lie group $G$, is given by a smooth map $\rho: G \times X \rightarrow X$, can be viewed as a group homomorphism 
$$\phi: G \rightarrow Diffeo(X),$$ given by $\phi(g)(x)=\rho(g,x)$ for all $g \in G$ and $x \in X$, where $Diffeo(X)$ is the group of all diffeomorphisms on $X$. Let $\mathfrak{g}$ be the Lie algebra of $G$ and $exp: \mathfrak{g} \rightarrow G$ be the exponential map.
Thus, the induced Lie algebra homomorphism $\tilde{\rho}: \mathfrak{g} \rightarrow \mathfrak{X}(X)$ is defined by $\xi \mapsto \tilde{\xi}$ where
\begin{align}\label{exp map}
	\tilde{\xi}(x):= \frac{d}{dt}(\rho(exp(t{\xi}), x))|_{t=0},
\end{align}
  for $x \in X,$ provides the associated fundamental vector field $\tilde{\xi}$ of $\xi \in \mathfrak{g}$.
It is well known that Lie derivative of a vector field $D_1$ along a vector field $D_2$ is given by the infinitesimal change of $D_1$ along the flow of $D_2$. Consider the one-parameter semi group of local diffeomorphisms of $X$ induced by the vector flow of $D_2$ by
$$\psi : (-\epsilon, \epsilon) \times X \rightarrow X~~\text{and}~~\psi_t(p):= \psi(t,p).$$
The Lie derivative of $D_1$ along $D_2$ at a point $x$ is defined as 
\begin{align*} \label{Lie derivative}
	(\mathfrak{L}_{D_2}D_1)(x)
	=&~ \underset{t \rightarrow 0}{lim} \frac{(\psi_{-t})_* D_1({\psi_t(x)})- D_1(x)}{t}\\
	=&~ \frac{d}{dt}|_{t=0}((\psi_{-t}))_* (D_1)_{\psi_t(x)}) 
\end{align*}
where $(\psi_t)_*$ is the pushforward along the diffeomorphism. See \cite{SM} for more details.

Similarly, for a Lie group action of $G$ on $X$, for each point $g \in G$ we have the diffeomorphism $\phi_g: X \rightarrow X$, induces an isomorphism $\tilde{\phi_g} : \mathfrak{X}(X) \rightarrow \mathfrak{X}(X)$ of $C^\infty(X)$-modules defined by the Lie derivative of a vector field along the flow given by $\phi_g$ is $($a conjugate action$)$ $$\tilde{\phi_g}(D_1)(x) = (\phi_{g})_*( (D_1)_{\phi_{g^{-1}}(x)}) =: g \cdot (D_1)_{\phi_{g^{-1}}(x)},$$
where $(\phi_g)_*$ is the pushforward of vector fields appears by the differential map of $\phi_g$.

But when go from the $C^\infty$ (smooth) real geometry to complex geometry and algebraic geometry, this type of information on the level of spaces consists of global sections does not work. There we need to deal with all local sections for which we use classical sheaf theoretic language accordingly. 
\subsection{Complex Lie group action on an analytic space} \label{G-analytic spaces}
Let $(X = \mathbb{C}^n, \mathcal{O}_X)$ be the standard complex manifold with its sheaf of holomorphic functions and $(Y = V(\mathcal{I}), \mathcal{O}_Y = \mathcal{O}_X/{\mathcal{I}})$ be an analytic space induced by a coherent ideal sheaf $\mathcal{I}\subset \mathcal{O}_X$. Denote the group of all biholomorphisms on $X$ by $biholo(X)$ and the group of all homeomorphisms on $Y$ by $homeo(Y)$.
\subsubsection{G-invariant analytic subspaces} \label{G-subspace} 
Let $G$ be a complex Lie group. Consider a holomorphic $G$-action on $X$, i.e. a group homomorphism 
$$\tilde{\phi}:G \rightarrow biholo(X).$$ The subspace $Y$ is said to be a $G$-invariant subspace of $X$ if the map $\tilde{\phi}$
 provides a group homomorphism $$\phi :G \rightarrow homeo(Y)~~\text{defined by}~~g \mapsto \phi_g:=\tilde{\phi}_g|_Y, ~~\text{for all}~~ g \in G.$$ We call $(Y,\mathcal{O}_Y)$ a $G$-analytic subspace of $X$ if in addition it satisfies the following criteria. 

If $U \subset X$ is open, then $V := Y\cap U$ is an open set of $Y$.
The $G$ action on $Y$, induces a $G$ action on $\mathcal{O}_Y$:  
$$G \times \mathcal{O}_Y(V) \rightarrow \mathcal{O}_Y(G\cdot V)$$
	as $(g, f)\mapsto g\cdot f$ where $ (g\cdot f)(y) = f(g^{-1}\cdot y)$,
compatible with restrictions. Note that $G \cdot V=\cup_{g \in G} (g \cdot V)$ is an open set in $Y$, since the homeomorphism $\phi_g:Y \rightarrow Y$ produces open sets $g \cdot V= \phi_g(V)$ for every $g \in G$ and arbitrary union of open sets is open.  This kind of group action on $Y$ canonically induces an $\mathcal{O}_Y - G$ module structure on $\mathcal{T}_Y=\mathcal{D}er_{\mathbb{C}_Y}(\mathcal{O}_Y)$ by the assignment $(g\cdot D)(f) := g\cdot (D(g^{-1}\cdot f))$, forms a $G$-action
	$$G\times Der_{\mathbb{C}}(\mathcal{O}_Y(V))\rightarrow Der_{\mathbb{C}}(\mathcal{O}_Y(G\cdot V)).$$
This action is compatible with the restriction morphism $res^{\mathcal{T}_Y}_{VW}$, for all $g\in G, D\in Der_{\mathbb{C}}(\mathcal{O}_Y(V))$ and for all open sets $W\subset V\subset Y$ (since for $g\in G, f\in \mathcal{O}_Y(g\cdot V)$ we get $g^{-1}\cdot f\in \mathcal{O}_Y(V)$ defined as $(g^{-1}\cdot f)(g^{-1}\cdot y) = f(y))$. Thus, there is a canonical homomorphism of sheaf of Lie algebras $$\tilde{\rho} :\mathfrak{g}_Y\rightarrow \mathcal{T}_Y$$ induced by the exponential map (\ref{exp map}) in the holomorphic context,
 where $\mathfrak{g}_Y$ and $\mathbb{C}_Y$ are the constant sheaves on $Y$ with stalks $\mathfrak{g}$ (the $\mathbb{C}$-Lie algebra of $G$) and $\mathbb{C}$ respectively.

\subsubsection{Induced action on $\Omega^1_Y$} We get a natural $\mathcal{O}_Y - G$ module structure on $\Omega^1_Y$ defined on local sections as 
$$(g\cdot \omega)(D) = g\cdot (\omega (g^{-1}\cdot D)),$$ for all $g\in G,~ D\in Der_{\mathbb{C}}(\mathcal{O}_Y(V)),~ \omega \in \Omega^1_Y(V)$ and for all open set $V\subset Y$.

	This process works for algebraic varieties instead of analytic spaces, for which we need to replace the standard complex Euclidean space $\mathbb{C}^n$ by the affine space $\mathbb{A}^n$ and sheaf of holomorphic functions replaced by sheaf of regular functions and biholomorphisms are replaced by biregular maps.

     But, the (equivalent) criteria of the $\mathcal{O}_Y - G$ module structure on $\mathcal{D}er_{\mathbb{C}}(\mathcal{O}_Y)$ does not holds here, since the exponential map is analytic but not a polynomial map (or regular map).

\subsection{Equivariant logarithmic de Rham cohomology} We begin by recalling the notion of equivariant de Rham cohomology and subsequently develop its analogue in the context of Lie algebroids.

\subsubsection{Equivariant (smooth) Lie algebroid cohomology} For a real $C^\infty$ manifold $X$ with a Lie group $G$ action,  the homotopy quotient or the homotopy orbit space $EG \times_G X$ is the orbit space $(EG \times X)/G$ of the diagonal $G$-action on $ EG \times X$ acts freely. the cohomology of the homotopy orbit space $EG \times_G X$ is known as the equivariant de Rham cohomology.
 It is given by the cohomology for the cochain complex formed by equivariant differential forms $\Omega_G^\bullet(X)$ together with the equivariant exterior derivative $$d_{G}: \Omega_G^\bullet(X) \rightarrow \Omega_G^{\bullet+1}(X)$$ (known as the Cartan model) given as follows (see \cite{GZ}):
 $\Omega_G^\bullet(X):=(Sym(\mathfrak{g}^*) \otimes_{\mathbb{R}} \Omega^\bullet(X))^G$ (i.e., forms by $G$-invariant elements of $Sym(\mathfrak{g}^*) \otimes_{\mathbb{R}} \Omega^\bullet(X)$), where $Sym(\mathfrak{g}^*)$ is the symmetric algebra of the dual Lie algebra $\mathfrak{g}^*$ of the Lie group $G$ and $\Omega^\bullet(X)$ is the graded algebra of differential forms on $X$,
	 and note that an equivariant differential form $\alpha$ is viewed as the polynomial map $\alpha: \mathfrak{g} \rightarrow \Omega^\bullet(X)$ defined by $\alpha(Ad(g)~D)=g~\alpha(D)$ for $g \in G$, $D \in \mathfrak{g}$, where $Ad: G \rightarrow Aut(\mathfrak{g})$ is the adjoint map. Thus, using the de Rham differential $d$  and the interior derivative $i_{\tilde{D}}$ by the fundamental vector field $\tilde{D} \in \mathfrak{X}(X)$ generated by $D \in \mathfrak{g}$, is
	$$d_{G}(\alpha)(D)=d(\alpha(D))-i_{\tilde{D}}(\alpha(D)).$$

Suppose $G$ is a compact, connected Lie group, and $X$ is a $G$-manifold. Then there is a canonical isomorphism,
$H_G^\bullet(X; \mathbb{R}):= H_{sing}^\bullet (EG \times_G X; \mathbb{R}) \cong H^\bullet((\Omega_G^\bullet(X), d_G))$.
Moreover, a $G$-manifold $X$ is said to be equivariantly formal (for example, $G$-acts trivially on $X$, see \cite{GZ}) if 
$$H_G^\bullet(X; \mathbb{R}) \cong (Sym(\mathfrak{g}^*))^G \otimes_{\mathbb{R}}  H^\bullet_{dR} (X),$$
as $(Sym(\mathfrak{g}^*))^G$-module, where $H^\bullet_{dR} (X)$ is the de Rham cohomology ring of the smooth manifold $X$.

\begin{Rem}
 By the Milnor's construction of the universal $G$-bundle $G \rightarrow EG \rightarrow BG:= EG/G$, the product space $EG \times X$ is homotopy equivalent to $X$ and the diagonal $G$-action on $EG \times X$ is a free action. 
Thus, the orbit space $X/G$ is homotopy equivalent to the homotopy quotient $EG\times_G X$ if the action of $G$ on $X$ is free.
Considering the associated $X$-bundle on $BG$, the Borel fibration $X \rightarrow EG \times_G X \rightarrow BG$, and using the associated spectral sequence, we get the above results.
\end{Rem}

 
 In \cite{BR}, a notion of equivariant cohomology for smooth Lie algebroids is introduced as follows.
 A Lie group $G$-action on a smooth Lie algebroid $\mathfrak{a}:L \rightarrow TX$ is given by a Lie algebra homomorphism $\rho: \mathfrak{g} \rightarrow \Gamma(L)$ such that the map $\Gamma(\mathfrak{a}) \circ \rho: \mathfrak{g} \rightarrow \mathfrak{X}(X)$ provides a $G$-action on $X$. Consider the graded vector space
 $$\mathcal{A}^\bullet=Sym(\mathfrak{g}^*) \otimes \Gamma (\wedge^\bullet L^*),$$
 together with the equivariant differential $d_{\mathfrak{g}}: \mathcal{A}^\bullet \rightarrow \mathcal{A}^{\bullet +1}$ by setting
 $$d_{\mathfrak{g}}(\mathcal{P} \otimes \omega)(\xi)=\mathcal{P}(\xi)(d_L(\omega)-i_{\rho(\xi)}(\omega)),$$
 where $\mathcal{P} \in Sym(\mathfrak{g}^*)$ and $\omega \in \Gamma (\wedge^\bullet L^*)$, $d_L$ is the Lie algebroid differential, $i_{\rho(\xi)}$ is the contraction operator with respect to ${\rho(\xi)}$ and both sides have evaluated on an element $\xi \in \mathfrak{g}$. Then the cohomology of the cochain complex $(\mathcal{A}_{\mathfrak{g}}^\bullet, d_{\mathfrak{g}})$ where $\mathcal{A}_{\mathfrak{g}}^\bullet:= {K}er~ d_{\mathfrak{g}}^2$ is denoted by $H_G^\bullet({L})$ and called the equivariant cohomology of the pair $({L}, \rho)$. This generalizes the classical case of $L= TX$. 
\subsubsection{Equivariant logarithmic de Rham complex and its hypercohomology} \label{equiv log de Rham} Here we extend the above ideas in the complex algebro-geometric settings to consider equivariant logarithmic de Rham cohomology. 

For a $G$-analytic space $Y$ of a Hermitian manifold $X$, we can canonically induce the cochain complex of sheaves of logarithmic differential forms (see (\ref{log differntial forms}))  $$(\Omega^\bullet_X(log~Y) \otimes_{\mathcal{O}_X} \mathcal{N}_{Y/X}, d)$$ (with poles along the divisor $Y$ and coefficient in the normal sheaf $\mathcal{N}_{Y/X}$ of $Y$ in $X$). If $Y$ is a (complex) submanifold, then we show this cochain complex of sheaves is $G$-equivariant, i.e. $d \circ (g \cdot ~) = (g \cdot ~) \circ d$ for each $g \in G$. 
In \cite{EE}, it has  been shown that for a module $(M, \nabla)$ over the Lie-Rinehart algebra $Der_{\mathbb{K}}(R)$, if the connection $\nabla$ is $G$-invariant, i.e. $g \cdot \nabla=\nabla$ for all $g \in G$ holds then $G$ acts on the Lie-Rinehart cohomology ring $H^\bullet(Der_{\mathbb{K}}(R), M)$, i.e. $g \cdot d= d \cdot g$ for any $g \in G$. Here, we consider the Lie-Rinehart algebra $\mathfrak{X}^T$ (usually denoted by $Der_{\mathbb{K}}(-log(I))$) as a local model and replace a module $(M, \nabla)$ by space of sections of a vector bundle $E$ with a covariant connection $\nabla$. The $G$-action for the covariant connection $\nabla$ on $E$ is defined as
$$(g \cdot \nabla)_D(s)=g \cdot \nabla_{g^{-1}\cdot D}(g^{-1}\cdot s),$$ for any $g \in G$ and for sections $D \in \mathcal{T}_X(U)$ and $s \in \Gamma(E|_U)$ over an open set $U$ of $X$.

\begin{Thm} \label{$G$-invariant}
	Any standard covariant connection on a $($smooth or holomorphic$)$ vector bundle over a $($real smooth or complex analytic$)$ $G$-manifold is $G$-invariant.
\end{Thm}
\begin{proof}
	Let $(X, \mathcal{O}_X)$ be a $($smooth or complex$)$ manifold with a $($complex$)$ Lie group $G$-action and $E \rightarrow X$ be  a vector bundle of rank $r$. Then the sheaf of sections $\Gamma_E$ is a locally free $\mathcal{O}_X$-module of rank $r$. Thus there is a local basis $\{s_1, \dots, s_r\}$ of $\Gamma({E})$ or $\Gamma(E|_U)=\langle\{s_1, \dots, s_r \}\rangle$ for some open set $U \subset X$  where $E|_U \rightarrow U$ is a trivial bundle. Then any standard covariant connection is locally defined as described in (\ref{covariant derivative locally}).
    For any $g \in G$, we need to show $(g \cdot \nabla)_D(s)=\nabla_D(s)$ holds, for all $D \in \mathcal{T}_X(U)$ and $s\in \Gamma(E|_U)$. Note that there exist $f_1, \dots, f_r$ in $\mathcal{O}_X(U)$ such that $s=\sum_{i=1}^{r}f_i s_i$.
Therefore, by definition we obtain
$$g^{-1} \cdot (\sum_{i=1}^{r}f_i s_i)= \sum_{i=1}^{r}(g^{-1}\cdot f_i)\cdot (g^{-1} \cdot s_i).$$ Thus, 
	\begin{align*}
		(g \cdot \nabla)_D(s)
		=&~g \cdot \nabla_{g^{-1}\cdot D}(\sum_{i=1}^{r}(g^{-1}\cdot f_i)~ (g^{-1} \cdot s_i))\\
		=&~ g \cdot (\sum_{i=1}^{r}((g^{-1} \cdot D)(g^{-1} \cdot f_i)) ~ (g^{-1} \cdot s_i)) \\
		=&~g \cdot (\sum_{i=1}^{r} (g^{-1} \cdot D(g(g^{-1} \cdot f_i))) ~ (g^{-1}\cdot s_i)) \\
		=&~g \cdot (\sum_{i=1}^{r} (g^{-1}\cdot D(f_i)) ~ (g^{-1}\cdot s_i))\\
		=&~g \cdot (\sum_{i=1}^{r} g^{-1} \cdot (D(f_i)~ s_i))\\
		=&~\sum_{i=1}^{r} D(f_i) \cdot s_i = \nabla_D(\sum_{i=1}^{r} f_i s_i)= \nabla_D(s).
	\end{align*}
\end{proof}
\begin{Cor}
	If a free divisor $Y$ is a $G$-analytic subspace of $X$, then the locally free $\mathcal{O}_X$-modules $\mathcal{T}_X(-log~ Y)$ and $\mathcal{N}_{Y/X}$ have canonical $G$-invariant connections.
\end{Cor}
\begin{Cor}
	The Lie algebroid cohomology $\mathbb{H}^\bullet(X,\Omega^\bullet_X(log~ Y)\otimes_{\mathcal{O}_X}\mathcal{N}_{Y/X})$ for a free divisor $Y$ with a $G$-analytic space structure $($see (\ref{log differntial forms})$)$, is $G$-equivariant.	
\end{Cor}
Using the notion of equivariant Lie algebrioid cohomology \cite{BR}, holomorphic Lie algebroid cohomology \cite{BRT} and above results, we define the equivariant Lie algebrioid complex for the Lie algebroid $\mathcal{T}_X(-log~Y)$ as
 $$(\mathcal{S}ym(\mathfrak{g}_X^*) \otimes_{\mathcal{O}_X} \Omega^\bullet_X(log~Y),~ \tilde{d}), $$  
and refer it by equivariant logarithmic de Rham complex and then consider its hypercohomology, an equivariant analogue of the logarithmic de Rham cohomology \cite{Ab}.

Furthermore, one possibility is extending the above cohomology to arbitrary (holomorphic) Lie algebroids, which will be referred as  \emph{equivariant $($holomorphic$)$ Lie algebroid cohomology}. This needs to define complex Lie group $G$ action on a holomorphic Lie algebroid and theory of holomorhic $G$-invariant connections.

\vspace{.5 cm}
\textbf{Acknowledgments.}
 I wholeheartedly express my deepest gratitude to my former PhD thesis supervisor Dr. Ashis Mandal for his immense support and guidance. Also, I would like to thanks Prof. Mainak Poddar for some of his valuable comments. The author also acknowledges support from the institutes IIT Kanpur, ISI Bangalore, and IISER Pune, as the work was partially conducted at each of these institutes to completion.

\vspace{.3 cm}
{\bf Abhishek Sarkar},\\
Department of Mathematics,\\
Indian Institute of Science Education and Research Pune \\
e-mail:  abhisheksarkar49@gmail.ac.in

\end{document}